\documentclass[10pt]{article}
\usepackage[latin2]{inputenc}
\usepackage{amsmath}
\usepackage[pdftex]{graphicx}
\usepackage{amsfonts}
\usepackage{amssymb}
\usepackage{algorithm}
\usepackage{algorithmic}
\usepackage[extra]{tipa}
\usepackage[left=1.5in]{geometry}
\usepackage{geometry}
\usepackage{wrapfig}
\usepackage{setspace}
\usepackage{epsfig}
\usepackage{color}
\usepackage{longtable}
\usepackage{listings}
\usepackage{tikz}
\usepackage{cancel}
\usepackage{float}
\usepackage{cite}
\usepackage{url}
\usepackage{hyperref}
\usepackage{amssymb, amsmath, epsfig, euscript, graphicx, xspace, amsthm, amscd, amsmath, mathtools}
\usepackage{float}
\restylefloat{table}
\usepackage{amsthm}
\usepackage{array,float}
\usepackage{setspace}
\usepackage{calligra}
\usepackage{enumitem}
\usepackage{tikz}
\usetikzlibrary{decorations.pathreplacing,calligraphy}
\usepackage{graphicx}
\usepackage{fancyref}
\usepackage{caption}
\usetikzlibrary{arrows,automata,positioning}
\usetikzlibrary{arrows}

\begin{document}

\newcommand{\Div}{\text{Div}}
\newcommand{\Pic}{\text{Pic}}
\newcommand{\Span}{\text{Span}}
\newcommand{\Edg}{\text{Edg}}
\newcommand{\ZZ}{\mathbb{Z}}
\newcommand{\QQ}{\mathbb{Q}}
\newcommand{\RR}{\mathbb{R}}
\newcommand{\CC}{\mathbb{C}}
\newcommand{\FF}{\mathbb{F}}
\newcommand{\NN}{\mathbb{N}}
\newcommand{\Cal}{\mathcal}
\newcommand{\mubar}{\overline{\mu}}
\newcommand{\bye}{\end{document}}
\newcommand{\aenum}{\begin{enumerate}[label=(\alph*)]}
\newcommand{\nenum}{\begin{enumerate}[label=(\arabic*)]}
\newcommand{\renum}{\begin{enumerate}[label=(\roman*)]}
\def\nor{\trianglelefteq}
\def\iso{\cong}
\def\Aut{\text{Aut}}
\def\QQ{\Bbb Q}
\def\ZZ{\Bbb Z}
\def\RR{\Bbb R}
\def\CC{\Bbb C}
\def\NN{\Bbb N}
\def\FF{\Bbb F}
\def\ord#1{| \, #1 \, |}
\def\gp#1{\langle \, #1 \, \rangle}
\def\Div{\text{Div}}
\def\Pic{\text{Pic}}
\def\vPic{\text{vPic}}
\def\Pr{\text{Pr}}
\def\RB{R\mathcal B}
\def\Gal{\text{Gal}}
\def\nor{\trianglelefteq}
\def\iso{\cong}
\def\Aut{\text{Aut}}
\def\QQ{\Bbb Q}
\def\ZZ{\Bbb Z}
\def\RR{\Bbb R}
\def\CC{\Bbb C}
\def\NN{\Bbb N}
\def\FF{\Bbb F}
\def\ord#1{| \, #1 \, |}
\def\gp#1{\langle \, #1 \, \rangle}
\def\Div{\text{Div}}
\def\Pic{\text{Pic}}
\def\vPic{\text{vPic}}
\def\Pr{\text{Pr}}
\def\RB{R\mathcal B}
\def\Gal{\text{Gal}}
\def\Char{\text{Char}}
\def\bigmid{\, \big | \,}
\renewcommand{\epsilon}{\varepsilon}
\newcommand{\DD}{\mathbb{D}}
\newtheorem{Theorem}{Theorem}
\newtheorem{Lemma}{Lemma}
\newtheorem{Corollary}{Corollary}
\newtheorem{proposition}{Proposition}
\newtheorem{definition}{Definition}
\theoremstyle{definition}
\newtheorem{exmp}{Example}

\begin{center}
\textbf{IWASAWA THEORY OF JACOBIANS OF GRAPHS}\\
\vspace{10px}
\textbf{SOPHIA R GONET}\\
\end{center}

\noindent \textbf{Abstract}
The Jacobian group (also known as the critical group or sandpile group)  is an important invariant of a finite, connected graph $X$;  it is a finite abelian group whose cardinality is equal to the number of spanning trees of $X$ (Kirchhoff's Matrix Tree Theorem). A specific type of covering graph, called a \textit{derived graph}, that is constructed from a \textit{voltage graph} with \textit{voltage group} $G$ is the object of interest in this paper.
Towers of derived graphs are studied by using aspects of classical Iwasawa Theory (from number theory). Formulas for the orders of the Sylow $p$-subgroups of Jacobians in an infinite voltage $p$-tower, for any prime $p$, are obtained in terms of classical $\mu$ and $\lambda$ invariants by using the decomposition of a finitely generated module over the Iwasawa Algebra.

\section{Introduction}

The Jacobian (or critical group, or sandpile group) is an algebraic invariant of a graph $X$ (in this paper the term graph will mean a simple graph with no loops or multiple edges, unless otherwise explicitly noted) which, for connected $X$, is a finite abelian group whose size is equal to the number of spanning trees of $X$ (this is well-known as the Matrix Tree Theorem). The study of Jacobians of graphs has a long history, and many applications, as described in \cite{19, 12, 3, 4, 2, 5, 24}.
Overall, there are relatively few graphs or families of graphs for which the Jacobian is exactly known: see \cite{1, 20, 8,18, 14, 9, 13}. In this paper we establish the ``asymptotic structure'' and orders of the Sylow $p$-subgroups of the Jacobians of certain covering graphs of a fixed base graph $X$,
namely those that belong to a cyclic voltage $p$-tower cover of $X$.

More specifically, we adapt to voltage towers of graphs the classical work of 
Iwasawa for $\ZZ_p$-extensions---infinite extensions $K_\infty$ of a number field $K$ with Galois group isomorphic to the additive $p$-adic integers, $\ZZ_p,$ for some prime $p$. 
By using the general theory of $\ZZ_p[[\Gamma]]$-modules, where $\Gamma=Gal(K_\infty/K)$, Iwasawa was able to unravel the structure of the inverse limit of the $p$-Sylow subgroups of the class groups of the finite extension fields in his towers. 
This enabled him to prove the following theorem, which can be found in \cite{28} and \cite{Iwa59}:
Let $K_\infty /K$ be a $\ZZ_p$-extension. Let $p^{e_m}$ be the exact power of $p$ dividing the order of the class group of $K_m,$ where $K_m$ is the fixed field of the subgroup $\Gamma^{p^m};$ then there exist nonnegative integers $\lambda, \mu$ and an integer $\nu$ such that $e_m=\mu p^m+\lambda m +\nu$ for all $m\geq m_0$ for some $m_0\geq 0$.\\

\noindent The Main Theorem of this paper is the analog in the graph theory setting:
\begin{Theorem}\label{main}
Let 
$$
X=X_0\leftarrow X_1\leftarrow X_2 \leftarrow \cdots \leftarrow X_m \leftarrow \cdots
$$
be a cyclic voltage $p$-tower (see Definition \ref{p-tower}), where all $X_m$ are connected. Let $\mathcal{J}_p(X_m)$ be the Sylow $p$-subgroup of the Jacobian of $X_m.$ Then there are nonnegative integers $\mu$ and $\lambda$ and an integer $\nu$ such that 
$$
|\mathcal{J}_p(X_m)|=p^{e_m}\qquad \text{where} \qquad e_m=\mu p^m +\lambda m+\nu
$$
for all $m\geq m_0$ for some $m_0\geq 0$.\\
\end{Theorem}

This theorem gives not only ``asymptotic'' order formulas for the $p$-Jacobians of the covering graphs $X_m$, but also 
their ``asymptotic'' invariant factor decompositions, which, in particular give conditions under which the $p$-ranks grow without bound (which is also analogous to the classical number theoretic results of Iwasawa). 

The ideas in \cite{17}, \cite{27} and \cite{30}
inspired the research that culminates in Section \ref{sec 4}. However, the work is independent, contemporaneous, and by quite different methods. 
%In contrast with Theorem \ref{main}, in \cite{27} Valli\`eres considers a Galois tower over a base graph $X$ that is allowed to be a multigraph; however, his results also require $X$ to be a {\it regular} multigraph. Using zeta-function and character-theoretic methods,  Valli\`eres obtains both {\it upper and lower asymptotic bounds} for the orders of the Sylow $p$-subgroups of the Jacobians in his towers, rather than an {\it exact} asymptotic (Iwasawa-type) formula.  [Note: Valli\`eres calls the towers ``abelian $p$-towers'' although they are the same as cyclic $p$-towers, since his definition also requires the Galois group of $X_m/X$ to be {\it cyclic} of order $p^m$.] So this work is complementary to \cite{27}, and uses different strategies and methodologies.

The terminology and theory of Jacobians, voltage graphs and their derived covering graphs --- including extending these 
results to infinite towers---is first summarized in Section \ref{sec 2}. In Section \ref{sec 4}, the theory of 
finitely generated modules over the Iwasawa algebra, $\ZZ_p[[\Gamma]]$ is summarized. 
Consequences of the latter results, that form the essential underpinning of the Main Theorem, are also established.
The Main Theorem is then proved using group-theoretic methods, that, in hindsight, illustrate 
how the decomposition theorem for Iwasawa modules plays the analogous role to the Smith Normal Form decomposition
that describes ordinary Jacobians.
The $p$-rank result mentioned above appears as a corollary to the Main Theorem.

This paper comprises the last part of the author's dissertation \cite{29}, which contains significantly more details, examples, and an array of additional theoretical and computational material on voltage graphs and their associated derived graphs.
We refer to it at points were its material expands on or expedites the development of this paper.

Added after refereeing: Just as this paper was submitted, Daniel Valli\`eres and Kevin McGown circulated a manuscript giving the 
generalization of Theorem \ref{main} to multigraphs, \cite{MVIII}.  
Their work---which is completely independent---uses ``analytic'' methods ($L$-series etc.),
and so provides a valuable complementary perspective on our result. It seems that the methods herein should also generalize to multigraphs, mutatis mutandis, 
since the main part of the proof, Section \ref{sec 4}, essentially only involves cokernels of Laplacians, and these are
well-defined for multigraphs.

\section{Preliminaries}
In Section \ref{2.1}, we define the Picard and Jacobian groups. Then in Section \ref{2.2} we define the Laplacian and reduced Laplacian. In Section \ref{3} we describe a specific type of covering graph, called a \textit{derived graph}, that arises from what is called a \textit{voltage graph}---where elements from a group (which may be finite or infinite) are assigned to the edges of a fixed base graph $X$. In Section \ref{3.4}, we give the definition of an \textit{intermediate covering graph}. We then state the important result: given a voltage graph with derived graph $Y$ such that $Y$ is connected, $Y/X$ is a normal (i.e., Galois) extension, and conversely, if $Y/X$ is a normal extension with Galois group $G$, then there exists a voltage assignment such that $(X,G,\alpha)$ is a voltage graph with derived graph $Y$. 
\label{sec 2}

\subsection{The Divisor, Picard, and Jacobian Groups}
\label{2.1}

For more details on this section, refer to \cite{7}.

\begin{definition}
A \textit{divisor} on a graph $X$ (possibly infinite) is an element of the free abelian group on the vertices $V=V(X)$:
$$
\Div(X) = \left\{ \sum_{v \in V(X)}a_v v \mid a_v \in \ZZ \right\}
$$
where each $\sum_{v \in V(X)} a_v v$ is a formal linear combination of the vertices of $X$ with
integer coefficients, where only finitely many $a_v$ are nonzero (in the case when $X$ is an infinite graph). The {\it degree} of a divisor is 
$$
\deg \left (\sum_{v \in V(X)} a_v v \right ) = \sum_{v \in V(X)} a_v  . 
$$
\end{definition}
\noindent When $V(X) = \{ v_1,\dots, v_n\}$, we may write the elements of $\Div(X)$ as 
$a_1 v_1 + a_2 v_2 + \cdots + a_n v_n$, where each $a_i \in \ZZ$,  and its degree is $a_1 + a_2 + \cdots + a_n$.
\\
\\
\noindent The \textit{degree map} $\deg:\Div(X)\to \ZZ,$ is a surjective group homomorphism with kernel equal to the subgroup of $\Div(X)$ of \textit{divisors of degree 0}, denoted as $\Div^0(X)$:
$$
\Div^0(X)=\{D\in \Div(X)\mid \deg D=0\}.
$$
Next let $X$ be a graph with vertices $\{v_1,\dots, v_n\}$.
For each fixed $v_i$ define the {\it principal divisor, $p_i$, based at $v_i$} by
$$
p_i=\deg(v_i)v_i-\sum_{j=1}^n \delta_{i,j}v_j
$$
where $\delta_{i,j}=1$ if $v_j$ is adjacent but not equal to $v_i$ and 0 otherwise 
(and here $\deg(v_i)$ is the valence of vertex $v_i$ in $X$).
Define {\it principal divisors} to be elements of the $\ZZ$-submodule of $\Div(X)$ spanned by the principal divisors based at the vertices:
$$
\Pr(X) = \text{Span}_{\ZZ}\{ p_i \mid 1 \le i \le n \}.
$$
Evidently $\Pr(X)$ is a submodule of $\Div^0(X)$.
From this we get the following groups.
\begin{definition}
The \textit{Picard group} of $X$ is the quotient group
$$
\Pic(X)=\Div(X)/\Pr(X),
$$
and the \textit{Jacobian group} of $X$ is the subgroup of $\Pic(X)$
$$
\mathcal{J}(X)=\Div^0(X)/\Pr(X).
$$
\end{definition}

\begin{Theorem}
If $X$ is connected, then $\mathcal{J}(X)$ is a finite abelian group.
\end{Theorem}

\subsection{The Laplacian and Reduced Laplacian}
\label{2.2}

%Firing moves may be expressed compactly via the graph Laplacian. 
\begin{definition}
Let $X$ be a graph with vertices $\{v_1,\dots, v_n\}.$ The graph \textit{Laplacian} $L=L_X$ is the $n\times n$ matrix given by
$$
L_{i,j}=\begin{cases}
\deg(v_i) &\text{ if }i=j\\
-1 &\text{ if } \text{ $v_i$ is adjacent to $v_j$}\\
0& \text{ if }\text{ $i\neq j$ and $v_i$ is not adjacent to $v_j$}
\end{cases}
$$
\end{definition}
%\noindent By Definition \ref{firing}, we see that the Laplacian matrix encodes all of the firing moves for $X$ since a firing move by vertex $v_j$ corresponds to subtracting the $j^{th}$ column of $L$ from a divisor.\\

\noindent The Laplacian is also the matrix representation of the following group homomorphism $\mathcal{L}$ defined as follows. 
$$
\mathcal{L}:\Div(X)\to \Div(X)\qquad \text{ where }\qquad \mathcal{L}(v_i)=p_i.
$$
When extended by $\ZZ$-linearity to all of $\Div(X)$, this is a $\ZZ$-linear homomorphism from $\Div(X)$ to itself, whose image is $\Pr(X)$, the group of principal divisors. From this, we get the following important fact:
$$
\Pic(X)=\Div(X)/\text{im}(\mathcal{L})=\text{coker}(\mathcal L).
$$
%where, by definition, the cokernel of a homomorphism $f:A\to B$ is $B/f(A).$ \\
%\\
%So it follows that the Picard group may be computed as the cokernel of the Laplacian. (For further details on this, see \cite{7}, Section 2.1.)\\
%\\
%We now relate this to the Jacobian and the \textit{reduced Laplacian}. Evidently, the matrix representation for $\mathcal{L}$ with respect to the basis of vertices is $L,$ the Laplacian matrix.\\
%\\
\noindent A \textit{reduced Laplacian} $\tilde{L}$ is the $(n-1)\times (n-1)$ integer matrix obtained by removing the row and column corresponding to any vertex $v$ from the Laplacian matrix $L$. So the Jacobian group can be computed as the cokernel of the reduced Laplacian matrix
$$
\mathcal{J}(X)\cong \ZZ^{n-1}/\text{im}(\widetilde{L})=\text{coker}(\widetilde{L}),
$$
where $\ZZ^{n-1}$ denotes the free $\ZZ$-module on the set $V(X)-\{v\}$ of rank $n-1.$

\subsection{Voltage Graphs}
\label{3}
We first give the definition of a general covering graph.
%We begin this section by defining a covering graph in general. 
%In Section \ref{3.1} we then go on to describe a specific type of covering graph, called a \textit{derived graph}, that arises from what is called a \textit{voltage graph}---where elements from a group (which may be finite or infinite) are assigned to the edges of a fixed base graph $X$.  In Section \ref{3.4}, we give the definition of an \textit{intermediate covering graph} $\widetilde{X}$. We then state several theorems that are important for Section \ref{sec 4}.  In Section \ref{3.5}, for $Y$ a derived graph (which may be infinite), we consider \textit{the group of divisors of $Y$}, \textit{the group of principal divisors of $Y$}, \textit{the Picard group of $Y$}, and \textit{the Laplacian of $Y$}, all as $\ZZ$ and $\ZZ[G]$ modules, where $G$ is any finite or infinite group.  \\

\begin{definition}
An undirected graph $Y$ is a \textit{covering} of an undirected graph $X$ if, after arbitrarily directing the edges of $X$, there is an assignment of directions to the edges of $Y$ and an onto graph homomorphism $\pi:Y\to X$ sending neighborhoods of $Y$ one-to-one onto neighborhoods of $X$ which preserve directions. We call such $\pi$ a covering map.
\end{definition}

\begin{definition}\label{dsheeted}
A \textit{$d$-sheeted} covering means every fiber contains exactly $d$ elements, i.e.,
 $$
 |\pi^{-1}(x)|=d \ \forall x\in V(X).
 $$
\end{definition}

	%\noindent We first start with an undirected, finite and connected graph $X$ and a group $G$ (which may be finite or infinite). We assume, for convenience that $X$ has no loops and no multiple edges; but the discussion easily generalizes to multigraphs. Next, we arbitrarily orient the edges of $X$. Then we label the forward-directed edges of $X$ with elements from the group $G$. (Note: the edges of $X$ do not have to labeled with different group elements. For instance, each edge could be assigned the group identity element.) These labels are referred to as the \textit{voltages} and the assignment itself is called the \textit{voltage assignment}. Furthermore, if we have a directed edge which goes from $u$ to $v$ labeled with the group element $\tau$, then label the edge which goes from $v$ to $u$ by the group element's inverse, namely $\tau^{-1}:$
\begin{definition} Let $X$ be a graph whose edges have been oriented, and let $G$ be a group (finite or infinite). For a fixed orientation of the edges of $X$, let $E(X)^+$ denote the set of forward-directed edges of $X$; and let $E(X)^-$ denote the same edges but each with the reverse orientation (so each undirected edge of $X$ becomes two edges in the disjoint union of $E(X)^+$ and $E(X)^-$). An (ordinary) \textit{voltage assignment} is a map
$$
\alpha: E(X)^+\cup E(X)^-\to G
$$
such that if $e_{i,j}\in E(X)^+$ and $\alpha(e_{i,j})=\alpha_{i,j}\in G,$ then $e_{j,i}\in E(X)^-$ and $\alpha(e_{j,i})=\alpha_{i,j}^{-1}$ (the inverse group element), where $e_{i,j}$ denotes the directed edge from $v_i$ to $v_j.$ \\
\\
The triple $(X, G, \alpha)$ is called an (ordinary) \textit {voltage graph}. The values of $\alpha$ are called the \textit{voltages} and $G$ is called the \textit{voltage group}.\\
\\
Note that a voltage assignment $\alpha$ is uniquely determined by its values on $E(X)^+$, so we will henceforth only specify $\alpha$ on the forward-directed edges of $X$.
\end{definition}

\noindent The vertices of $X$ are labeled as $v_1,\dots,v_n$.
This imposes a natural lexicographic orientation on $X$, namely whenever there is an edge between $v_i$ and $v_j$,
orient the edge $v_i \rightarrow v_j$ if $i < j$ (called the \textit{standard orientation}). Note that results on derived graphs do not depend on the choice of orientation by \cite{29}, so without further mention, we adopt the standard orientation.\\
\\
 \noindent Any such voltage assignment can be codified by its $n \times n$ {\it voltage adjacency matrix}
$$
A_\alpha = \begin{pmatrix} A_{i,j} \end{pmatrix}
$$
with entries $A_{i,j}\in \ZZ[G]$ such that $A_{i,j}=0$ if $i=j$ or there is no edge between $v_i$ and $v_j$, and $A_{i,j}=\alpha_{i,j}$ otherwise. (Note that the {voltage adjacency matrix} is also defined in \cite{15} as Definition 2.16.)
%The entries of $A_\alpha$ are from the integral group ring $\ZZ [G]$, and $A_\alpha$
%is ``inverse symmetric'' in the sense that its transpose, $A_\alpha^t$, is the same as $A_\alpha$, but with all group elements in nonzero entries inverted. 
%(As usual, we identify the integer $k$ with the group ring element $k \tau_0$, where, for formality, $\tau_0$ is the identity of $G$.)
%Changing the orientation of $X$ simply results in a voltage assignment where some $\alpha_{i,j}$ are inverted. \\

\noindent The purpose of assigning voltages to the graph $X$, called \textit{the base graph}, is to obtain an object called \textit{the derived graph}, called $Y$ here. To get the vertices of $Y$, make $d=|G|$ copies of each vertex $x\in V(X)$ labeling them as $x_{\tau_0}, x_{\tau_1}, x_{\tau_2},\dots,x_{\tau_{d-1}}$ where $G=\{\tau_0,\tau_1, \tau_2,\dots, \tau_{d-1}\}$ has order $d$ (and the same formal construction works even if $|G|$ is uncountable). So there are $|G|\cdot |V(X)|$ vertices in $Y$. 
Now create the edges of $Y$ by the following rule:
whenever there is an edge from $v_i$ to $v_j$ in the base graph $X$ with assigned voltage $\alpha_{i,j}$, 
create edges that go from $v_{i,g}$ to $v_{j, g \alpha_{i,j}}$ in $Y$, for every $g \in G$,
where $g \alpha_{i,j}$ is the group-product of these two group elements in $G$. 
If $|G|=d,$ then $\pi:Y\to X$ is a $d$-sheeted covering map (where again, $d$ may be any infinite cardinal too). Note that the degree (valence) of each vertex $v_\tau$ of $Y$ is the same as the degree of $v=\pi(v_\tau)$ in $X$. Also, since our base graph $X$ has no loops (i.e., $i \ne j$ here), no two vertices in the same fiber of $\pi$ are adjacent in~$Y$.\\
\\
Many examples as well as computational ways of constructing the ordinary adjacency matrix of $Y$ from the voltage adjacency matrix of $X$ by tensoring with matrices for the regular representation of $G$ appear in \cite{29}.

\subsection{Galois Theory of Covering Graphs and Voltage Graphs}
\label{3.4}

\noindent Refer to \cite{26} for Galois theory of (finite) Galois covers. For proof of the theorems presented below, see \cite{29}.

\begin{definition}\label{inter}
Suppose $Y$ is a covering of $X$ with projection map $\pi.$ A graph $\widetilde{X}$ is an \textit{intermediate covering} to $Y/X$ if $Y/\widetilde{X}$ is a covering, $\widetilde{X}/X$ is a covering and the projection maps $\pi_1:\widetilde{X}\to X$ and $\pi_2:Y\to \widetilde{X}$ have the property that $\pi=\pi_1\circ \pi_2$. 
If $Y/X$ is a $d$-sheeted covering with projection map ${\pi}:Y\to X,$ then it is \textit{normal} or \textit{Galois} if there are exactly $d$ graph automorphisms ${\sigma}:{Y}\to Y$ such that ${\pi}\circ \sigma={\pi}.$ The \textit{Galois group} is $G=Gal(Y/X)=\{\sigma:Y\to Y \ | \ \pi\circ\sigma=\pi\}.$\\
\end{definition}
\begin{Theorem}\label{FTGTCor}
Suppose $Y/X$ is a normal covering with Galois group $G$ and $\widetilde{X}$ an intermediate covering corresponding to the subgroup $H$ of $G$. Then $\widetilde{X}$ itself is a normal covering of $X$ if and only if $H$ is a normal subgroup of $G$, in which case $Gal(\widetilde{X}/X)\cong G/H.$
\end{Theorem}

\noindent Now we put this in terms of voltage graphs.

\begin{Theorem}\label{theorem5}
Let $(X,G,\alpha)$ be a voltage graph with $Y$ the derived graph. If $Y$ is connected, then $Y/X$ is a normal cover with $Gal(Y/X)\cong G$. Conversely, given a normal (Galois) cover $Y/X$, with $G=Gal(Y/X),$ then $Y/X$ is a voltage cover with the voltage group equal to $Gal(Y/X).$
\end{Theorem}

\noindent For any Galois cover $\pi : Y \rightarrow X$ of a connected base graph $X$, the graph
$Y$ is necessarily connected 
except in the case where $G = \Gal(Y/X)$ is the cyclic group of order 2 and $Y$ is two disjoint isomorphic 
copies of $X$ interchanged by $G$.
\begin{Theorem}\label{voltage}
Let $(X,G,\alpha)$ be a voltage graph with derived graph $Y$ such that $Y$ is connected. If $\widetilde{X}$ is an intermediate cover of $Y/X$ corresponding to the normal subgroup $H$ of $G$, then $\widetilde{X}/X$ is a voltage graph, whose voltage adjacency matrix is the voltage adjacency matrix of $Y/X$, but with nonzero entries reduced modulo $H$ (thus has entries in $\ZZ[G/H]).$
\end{Theorem}
%\begin{proof}
%Let $\widetilde{X}$ be an intermediate cover of $Y/X$ that corresponds to the normal subgroup $H$ of $G$. So we have that $Y/\widetilde{X}$ and $\widetilde{X}/X$ are normal coverings with projection maps $\pi_1:\widetilde{X}\to X$ and $\pi_2:Y\to \widetilde{X}$ such that $\pi_2\circ \pi_1=\pi.$ \\
%\\
%Let $A_{\alpha,Y}$ be the voltage adjacency matrix for $X$ corresponding to the derived graph $Y$. So $A_{\alpha,Y}$ has entries in $\ZZ[G].$ Now reduce the entries of $A_{\alpha,Y}$ modulo $H$, and denote this matrix $A_{\alpha_1}$, which has entries in $\ZZ[G/H].$ This encodes the voltage graph $(X,G/H,\alpha_1)$, where $\alpha_1: E(X)^+\to G/H$. Since $Gal(\widetilde{X}/X)\cong G/H$ by Theorem \ref{FTGTCor}, it follows that the derived graph of $(X,G/H,\tilde{\alpha})$ is $\widetilde{X}$ by the Fundamental Theorem of Galois Theory for Graphs.
%\end{proof}

\section{Towers of Voltage Graphs and Iwasawa Theory}
\label{sec 4}

We begin Section \ref{4.1} by defining a \textit{cyclic $p$-tower} of graphs. We then extend this definition to a cyclic \textit{voltage} $p$-tower of graphs by using Theorems \ref{theorem5} and \ref{voltage} from Section \ref{3.4}. From this, we get a ``universal cover'' of the tower by an infinite derived graph that we call $X_{p^\infty}; $ it is the derived graph obtained from the voltage graph $(X,\ZZ_p,\alpha)$, where the voltage group is the additive $p$-adic integers and the voltage assignment $\alpha$ is determined by the cyclic voltage $p$-tower. We call $X_{p^\infty}$ the \textit{completion of the tower}. In Section \ref{4.3}, we present important definitions and results pertaining to $\Lambda$-modules. Then in Subsection \ref{4.3.1}, we specify the $\Lambda$-modules be finitely generated. In Subsection \ref{4.3.2} we construct a finitely generated torsion $\Lambda$-module, which we call $\Pic_\Lambda.$ Finally in Section \ref{4.4} we prove the main theorem (Theorem \ref{main}) of this paper.

\subsection{$p$-Tower Covering Graphs}
\label{4.1}

We begin by defining a cyclic $p$-tower of graphs. We assume $p$ is a fixed prime.\\
\begin{definition}\label{tower}
A \textit{cyclic $p$-tower} of graphs above a base graph $X$ is a sequence of covering graphs
$$
X=X_0\leftarrow X_1\leftarrow X_2 \leftarrow \cdots \leftarrow X_m \leftarrow \cdots
$$
such that for $m\geq 0,$ the cover $X_{m}/X$ is normal with $Gal(X_m/X)\cong\ZZ/p^m\ZZ.$
\end{definition}
\noindent Note that for $m\geq 0,$ this implies that the cover $X_{m+1}/X_m$ is normal with $Gal(X_{m+1}/X_m)\cong\ZZ/p\ZZ$ by the Fundamental Theorem of Galois Theory, along with the Third Isomorphism Theorem from \cite{11}.\\

\noindent Now we specialize Definition \ref{tower} to voltage graphs. 
\begin{definition}\label{p-tower}
A \textit{cyclic voltage $p$-tower} of graphs above a base graph $X$ is a sequence of derived graphs
$$
X=X_0\leftarrow X_1\leftarrow X_2 \leftarrow \cdots \leftarrow X_m \leftarrow \cdots
$$
such that
for $m\geq 0,$ $X_{m}/X$ is a derived graph with $Gal(X_m/X)\cong\ZZ/p^m\ZZ.$
\end{definition}
\noindent Theorems \ref{theorem5} and \ref{voltage} in Section \ref{3.4} extend to towers, 
so we may choose notation that describes the vertices, edges, voltage assignments and Galois actions on
the graphs $X_m$ in compatible ways that are determined by the covering maps. In short, each $X_m$ is a derived cover of $X$ defined by a voltage adjacency matrix of fixed degree $n$ 
whose $i,j$ entries are $\alpha_{i,j}^m$ in the group ring $\ZZ[\ZZ/p^m \ZZ]$; and for each fixed $i,j$ the sequence of such entries has a limit as $m \rightarrow \infty$ in the $p$-adic integral group ring. The $n \times n$ matrix whose entries 
are these limits forms a voltage adjacency matrix for the ``completion graph'' that we now describe.
This construction is straightforward, and the precise details are given explicitly in Section~5.1 of \cite{29}. 
This leads to a ``universal cover''
of the tower, by an infinite derived graph that we call~$X_{p^\infty}$.\

\begin{definition}\label{Xpinfty}
Given a cyclic voltage $p$-tower as in Definition \ref{p-tower}, with each $X_m$ the derived graph for the voltage assignment
$\alpha_m : E(X)^+ \rightarrow \ZZ/p^m \ZZ$, let $X_{p^\infty}$ be the derived graph obtained from
the voltage graph $(X, \ZZ_p,\alpha)$, where the voltage group is the (additive) $p$-adic integers and voltage assignment 
$\alpha$ is determined by the tower.
We call $X_{p^\infty}$ the \textit{completion of the tower}.
\end{definition}

\begin{Theorem}\label{towercompletion}
Let $X = X_0 \leftarrow X_1 \leftarrow X_2 \leftarrow \cdots \leftarrow X_m \leftarrow \cdots$ be a  
cyclic voltage $p$-tower, let $X_{p^\infty}$ be the completion of the tower, and for each $m \ge 0$ let $\overline X_m$
be the associated intermediate graphs.
Then for all $m \ge 0$ there are graph isomorphisms $\overline X_m \rightarrow X_m$, depicted as the horizontal maps in 
Figure \ref{figure 5.2}, such that all the maps in that figure commute and commute with the action of $\ZZ_p$ as automorphisms of each graph.

\end{Theorem}

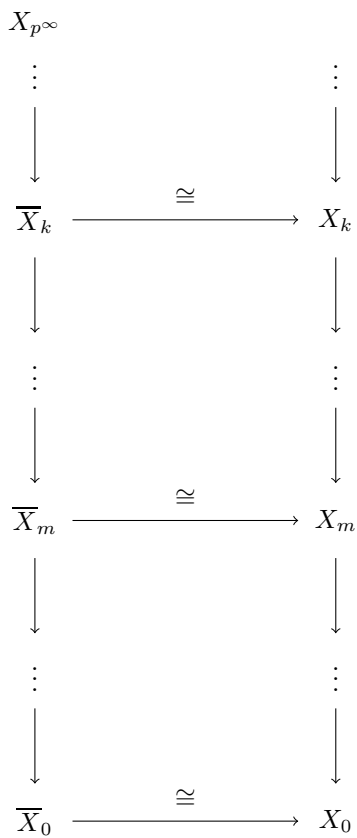
\begin{figure}[H]
\centering
\begin{tikzpicture}
%% vertices
%% vertex labels
\node at (0,0) {$\overline{X}_0$};
\node at (0,4) {$\overline{X}_m$};
\node at (0,8) {$\overline{X}_k$};
\node at (2,8.3) {$\iso$};
\node at (2,4.3) {$\iso$};
\node at (2,.3) {$\iso$};
\node at (0,10.6) {${X}_{p^\infty}$};
\node at (0,2) {$\vdots$};
\node at (0,6) {$\vdots$};
\node at (0,10) {$\vdots$};
\node at (4,0) {$X_0$};
\node at (4,4) {${X}_m$};
\node at (4,8) {${X}_k$};
\node at (4,2) {$\vdots$};
\node at (4,6) {$\vdots$};
\node at (4,10) {$\vdots$};
%%% edges
\draw[->] (0,1.5) -- (0,.5);
\draw[->] (0,3.5) -- (0,2.5);
\draw[->] (0,5.5) -- (0,4.5);
\draw[->] (0,7.5) -- (0,6.5);
\draw[->] (0,9.5) -- (0,8.5);
\draw[->] (4,1.5) -- (4,.5);
\draw[->] (4,3.5) -- (4,2.5);
\draw[->] (4,5.5) -- (4,4.5);
\draw[->] (4,7.5) -- (4,6.5);
\draw[->] (4,9.5) -- (4,8.5);
\draw[->] (.5,8) -- (3.5,8);
\draw[->] (.5,4) -- (3.5,4);
\draw[->] (.5,0) -- (3.5,0);
\end{tikzpicture}
\caption{Cyclic voltage $p$-tower with completion $X_{p^\infty}$}
\label{figure 5.2}
\end{figure}

\noindent We will henceforth write the voltage group $\ZZ_p$ as the 
{\it multiplicative} profinite group $\Gamma$, where, as a profinite group, it is cyclic:
it is the closure of an infinite (multiplicative) cyclic group $\gp \gamma$ under the $p$-adic metric topology, 
for some $\gamma$. \\
\\
%\noindent If we followed the ``Iwasawa program'' directly, we would try to show that the Sylow $p$-subgroup of $\mathcal{J_\infty}$ is a finitely generated torsion $\Lambda$-module. However, in the next few sections we take a somewhat different tack. Since our towers are voltage towers, we can use divisors to create modules over $\ZZ_p[\Gamma]$ and then over $\Lambda$ to construct a finitely generated torsion $\Lambda$-module that maps onto each finite Jacobian, which we will do in Subsection \ref{4.3.2}. In this way, we can ultimately prove our Main Theorem, Theorem \ref{main}, on the growth of the orders of these finite Jacobians. So we achieve the same goal as the Iwasawa program via this~route.\\
%\\
%For the moment we leave aside the question as to whether our Jacobian $\Lambda$-module, constructed via divisors actually coincides with the inverse limit of the Sylow $p$-subgroups of the finite Jacobians, $\mathcal{J}_p(X_{p^\infty})$, since there may be unresolved issues of inverse limits and tensor products.\\
%\\
%\noindent In the next section, we give important results and definitions pertaining to $\Lambda$-modules. Many of the modules over $\Lambda$ that occur in Iwasawa theory are finitely generated torsion modules. Thus, in Subsection \ref{4.3.1}, the $\Lambda$-modules that we consider will be finitely generated. 

\subsection{Iwasawa Modules}
\label{4.3}
Fundamental to Iwasawa's development of $p$-class groups in $\ZZ_p$-towers of number fields was his study of certain finitely generated modules over the $\ZZ_p$-algebra $\Lambda = \ZZ_p[[\Gamma]]$. Here $\Lambda$ is the compactification of $\ZZ_p[\Gamma]$, under the profinite topology defined by the open subgroups $\Gamma^{p^m}$. In this section we list and establish some properties of $\Lambda$, as well as finitely generated modules over it, that are the underpinnings of our main theorem.
The following two useful theorems about $\Lambda$ can be found in \cite{28}. 
\begin{Theorem}\label{25}
For the indeterminate $T$, $\ZZ_p[[\Gamma]]\cong \ZZ_p[[T]]$ with the isomorphism being induced by $\gamma\mapsto T+1.$ 
\end{Theorem}

\begin{Theorem}\label{noeth}
$\Lambda=\ZZ_p[[\Gamma]]$ is a Noetherian local ring.
\end{Theorem}

%\noindent We now define what it means for two $\Lambda$-modules to be \textit{pseudo-isomorphic}. We then define when a nonconstant polynomial is \textit{distinguished}. This leads the way to results (in particular, Proposition \ref{proposition9} and Corollary \ref{corollary6}) that will be used in proving the claims in Section \ref{4.4}.

\begin{definition}
Two $\Lambda$-modules $M$ and $M'$ are said to be \textit{pseudo-isomorphic}, written
$$
M\sim M'
$$
if there is a homomorphism $M\to M'$ with finite kernel and co-kernel.
\end{definition}
\begin{definition}
A nonconstant polynomial $P(T)\in \Lambda$ 
$$
P(T)=T^n+a_{n-1}T^{n-1}+\cdots+a_0
$$
is called \textit{distinguished} if $p\mid a_i$ for all $0\leq i\leq n-1.$
\end{definition}

\noindent The following proposition can be be obtained immediately from Proposition 7.2 and Lemma 7.5 of \cite{28}.
\begin{proposition}\label{proposition8}
Let $F(T)$ be a distinguished polynomial in $\ZZ_p[T]$.  Then
$$
\ZZ_p[T]/F(T)\ZZ_p[T] \iso \ZZ_p[[T]]/F(T)\ZZ_p[[T]],
$$
where the isomorphism is as $\ZZ_p[T]$-modules.
The isomorphism is the natural one, namely for $r \in \ZZ_p[T]$, the coset $r + F(T)\ZZ_p[T]$ maps to 
$r + F(T)\ZZ_p[[T]]$.
\end{proposition}

\bigskip

\noindent Let $\Lambda = \ZZ_p[[\Gamma]]$ and
fix a topological generator $\gamma$ for $\Gamma$;
so by Theorem \ref{25} the map $\gamma \mapsto T + 1$ extends to an isomorphism from $\Lambda$ to $\ZZ_p[[T]]$.
For $m \ge 0$ let $\omega_m = \gamma^{p^m} - 1$.

\begin{Lemma}\label{lemma1}
For $m \ge 0$ $\omega_m$ maps to a distinguished polynomial in $\ZZ_p[T]$.
\end{Lemma}

\begin{proof}
By definition, $\omega_m$ maps to $(T+1)^{p^m} - 1$.  
Thus
$$
(T+1)^{p^m} - 1 \equiv (T^{p^m} + 1^{p^m}) - 1 \equiv T^{p^m} \pmod{p\ZZ_p[T]},
$$
which establishes claim. 
\end{proof}
%\noindent Note: this works for $\omega_0 = \gamma - 1$ too.

\bigskip

\medskip

\noindent Let $R=\ZZ_p[\Gamma].$ Fix $m\geq 0.$

\begin{definition}\label{Omega}
Let $D$ be any $\Lambda$-module and let $B$ be any subset of $D$. For every $m\geq 0,$ define 
$$
\Omega^D_m(B) = B\cap \omega_m D .
$$
Define $R_m = R / \Omega^\Lambda_m({R})=R/R\cap \omega_m \Lambda$.
\end{definition}
%\noindent Note: with $\omega_m$ in place of $\eta_m,$ we define $\mho(B)=B\cap \eta_m A$.
\noindent In the special case when $B$ is an $R$-submodule of $D$ (where $D$ is considered as an $R$-module),
we have that $\Omega^D_m(B)$ is an $R$-submodule of $B$ containing $\omega_m B$.\\
\\
The sets $\Omega^D_m(B)$ define relatively open subsets of $B$ in the ``$\omega$-adic topology'' on $D$.
They obey the appropriate transitive property:
If $B$ and $C$ are subset of $D$ with $C \subseteq B$, then
$$
\Omega^D_m(B) \cap C = \Omega^D_m({C}).
$$
It is not true in general that $\omega_m B \cap C = \omega_m C$ however.

\bigskip

\begin{proposition}\label{proposition9}
Let $D$ be a $\Lambda$-module, 
let $A$ be any $\Lambda$-submodule of $D$ and let $B$ be any $R$-submodule of $A$, where $A$ is considered as an $R$-module.
Then the map
$$
\phi : B / \Omega^D_m(B) \longrightarrow A/\Omega^D_m(A) \qquad\text{by}\qquad
\phi(x + \Omega^D_m(B)) = x + \Omega^D_m(A)
$$
is a well-defined and injective $R$-module homomorphism. 
If $B$ contains a set of $\Lambda$-module generators for $A$, then $\phi$ is an isomorphism 
and $A = B + \Omega^D_m(A)$; and if additionally 
$\omega_m D \subseteq A$ then $A = B + \omega_m D$.
\end{proposition}

\begin{proof}
We first simplify notation by denoting $\Omega^D_m({C})$
by just $\Omega_m({C})$ for every subset $C$ of $D$ throughout the proof. The map $B \rightarrow A/\Omega_m(A)$ by $x \mapsto x + \Omega_m(A)$ is a well-defined $R$-module homomorphisms, and 
since $\Omega_m(B) \subseteq \Omega_m(A)$, its 
kernel clearly contains $\Omega_m(B)$.  This map therefore factors through $B / \Omega_m(B)$, giving the homomorphism $\phi$.
We also have
\begin{align*}
\ker \phi &= \{ x + \Omega_m(B) \mid x \in B \text{ and } x + \Omega_m(A) = 0 + \Omega_m(A)  \} \\
&= \{ x + \Omega_m(B) \mid x \in B \text{ and } x \in \Omega_m(A) \} \\
&= (B \cap \Omega_m(A))/\Omega_m(B) \\
&= (B \cap (A \cap \omega_m D)) / \Omega_m(B) \\
&= (B \cap \omega_m D) / \Omega_m(B) 
= \Omega_m(B)/\Omega_m(B) = 1,
\end{align*}
so $\phi$ is injective.
It remains to show 
if $B$ contains a set of $\Lambda$-module generators for $A$, then $\phi$ is surjective.
Assuming this hypothesis, every $y \in A$ can be written as 
$$
y = \alpha_1 b_1 + \cdots + \alpha_n b_n, \qquad\text{for some } \alpha_1,\dots,\alpha_n \in \Lambda \text{ and } b_1,\dots,b_n \in B.
$$
By Proposition~\ref{proposition8} and Lemma~\ref{lemma1}, for each $\alpha_i$ there is some $r_i \in \ZZ_p[\gamma] \subseteq R$ 
such that $\alpha_i - r_i \in \omega_m \Lambda$.
Let $y' = r_1 b_1 + \cdots + r_n b_n \in B$.
By construction,
$$
y-y' = (\alpha_1 - r_1) b_1 + \cdots + (\alpha_n - r_n) b_n \in A \cap \omega_m D = \Omega_m(A).
$$
Thus $\phi(y' + \Omega_m(B)) = y' + \Omega_m(A) = y + \Omega_m(A)$, and so $\phi$ is surjective, hence an isomorphism. Also, surjectivity of $\phi$ implies that $A = B + \Omega_m(A)$.
If $\omega_m D \subseteq A$, then $\Omega_m(A) = \omega_m D$, so the last assertion holds too.

\end{proof}

\begin{Corollary}\label{corollary6}
For $R=\ZZ_p[\Gamma]$, 
we have that $\phi$ induces an $R$-module isomorphism
$$
R_m \iso \ZZ_p[\Gamma_m].
$$
\end{Corollary}

\begin{proof}
As in Definition \ref{Omega}, we have $R_m=R/\Omega_m({R})=R/(R\cap\omega_m\Lambda)$. Now  $\ZZ_p[T]$ corresponds to $\ZZ_p[\gamma]$ in the isomorphism between $\ZZ_p[[T]]$ and $\Lambda$, where $\gamma$ is a fixed topological generator for $\Gamma$. So we have 
\begin{align*}
R_m &= R/(R \cap \omega_m \Lambda)&& \text{(by definition)}\\
&\cong \Lambda / \omega_m \Lambda && \text{(by Proposition \ref{proposition9})}\\
&\cong \ZZ_p[\gamma]/(\omega_m)&& \text{(by Proposition \ref{proposition8})}\\
&\cong \ZZ_p[\gamma]/(\gamma^{p^m}-1) \\
&\cong \ZZ_p[\Gamma_m],
\end{align*}
where the last isomorphism follows since $\ZZ_p[\gamma]/(\gamma^{p^m}-1)$ is isomorphic to the group ring of the cyclic group $\ZZ/p^m \ZZ \iso \Gamma_m$. Hence $R_m\cong \ZZ_p[\Gamma_m].$

\end{proof}

\noindent We now record an elementary lemma, which will be used in proving Claim 6(3) in Section \ref{4.4}.
\begin{Lemma}\label{ZZp}
If $A \iso \ZZ_p$ as a $\ZZ_p$-module and $B$ is a $\ZZ_p$-submodule of $A$, then either $B = 0$ or $A/B$ is finite.
\end{Lemma}
\begin{proof} 
By hypothesis $A$ is isomorphic to the ring $\ZZ_p$ considered as a module over itself,
so its submodules are ideals.
If $B \ne 0$, then $B  = p^k A$, for some $k \ge 0$, and so $A/B \iso \ZZ_p/p^k \ZZ_p \iso \ZZ/p^k \ZZ$, which is finite.
\end{proof}

\noindent In the next subsection, we now specify our $\Lambda$-modules to be finitely generated. Theorem \ref{growth} will be used in Section \ref{4.4} to ultimately prove Theorem \ref{main}.

\subsubsection{Finitely Generated $\Lambda$-Modules}
\label{4.3.1}

\noindent The following theorem can be found in \cite{28} as Theorem 13.12. 
\begin{Theorem}\label{structure2}[Structure Theorem for Iwasawa modules]
For any finitely generated $\Lambda$-module $M,$ we get the following pseudo-isomorphism:
$$
M\sim \Lambda^r \oplus \left(\oplus_{i=1}^s \Lambda /(p^{k_i}) \right) \oplus \left(\oplus_{j=1}^{t} \Lambda/(g_j(T)^{m_j})\right)
$$
where $r=\text{rank}(M),$ $s, t,k_i$ and $m_j\in \ZZ$ and $g_i \in \ZZ_p[T]$ is monic, distinguished and irreducible. This decomposition is uniquely determined by $M$. If $M$ is a torsion module, then $r=0.$
\end{Theorem}
\noindent The growth formula for the orders of the finite Jacobians in the conclusion of the main result, Theorem \ref{main},
ultimately comes from the orders of certain finite quotients of the cyclic factors in the Iwasawa Structure Theorem decomposition of
a finitely generated, torsion Iwasawa module that we shall construct shortly. 
The general structure of finite quotients of cyclic $\Lambda$-modules is described in \cite{28}, Section 13.3.

\begin{definition}
As in the notation of Theorem \ref{structure2}, we define the \textit{Iwasawa invariants} of $M$ by
$$
\mu=\sum_{i=1}^s k_i \qquad \text{and } \qquad \lambda=\sum_j m_j \deg g_j.
$$
\end{definition}
\begin{definition}\label{characteristic}
Let $M$ be any finitely generated torsion $\Lambda$-module with $p^{k_i}$ and $g_j^{m_j}$, as in Theorem \ref{structure2}.
The \textit{characteristic polynomial} of $M$, denoted by $\Char(M)$, is the product:
$$
\Char(M) = p^{k_1 + \dots + k_s} g_1^{m_1} \cdots g_t^{m_t},
$$
where $\Char(M) = 1$ if $M$ is finite.
\end{definition}

\noindent We record some basic facts about finitely generated torsion $\Lambda$-modules.
These may be found in Section 1.1 of Bence Forr\'as Master's Thesis ``Iwasawa Theory,'' \cite{For20}. Part (3) may also be found in \cite{CS}. 

\begin{proposition}\label{forras}
Let $P$ be any finitely generated torsion $\Lambda$-module.
\begin{enumerate}
\item
The relation ``pseudo-isomorphism'' is an equivalence relation on any set of finitely generated torsion $\Lambda$-modules.

\item
For any $\Lambda$-module $M$, the characteristic polynomial is an invariant of the pseudo-isomorphism equivalence class of $M$.

\item
If $M$ is a submodule of $P$, then $\Char(P) = \Char(M) \Char(P/M)$. In particular, $\Char(M) \bigmid \Char(P)$.
\end{enumerate}

\end{proposition}
\noindent The following Theorem can be found in Romyar Sharifi's online notes ``Iwasawa Theory'' as Theorem 2.4.7, \cite{Shar}. In it, $\mathcal{O}$ is a valuation ring of a $p$-adic field. We simplify his statement by taking $\mathcal{O}=\ZZ_p.$ This can also be obtained from Proposition 13.19 and Lemma 13.21 in \cite{28}.
\begin{Theorem}\label{Sharifi} 
Let $M$ be a finitely generated, torsion, $\Lambda$-module, and let $n_0\geq 0$ be such that $\Char(M)$ and $\omega_{m,n_0} = \omega_m / \omega_{n_0}$ are relatively prime for all $m\geq n_{0}.$ Set $\lambda(M) =\lambda$ and $\mu(M)=\mu.$ Then there exists an integer $\nu$ such that
$$
|M/\omega_{m,n_0}M| = q^{e_m} \qquad \text{where} \qquad e_m = \mu p^m +\lambda m +\nu
$$
for all sufficiently large $m\geq 0.$
\end{Theorem}
\noindent This theorem is used to prove Theorem \ref{growth}. First we present two lemmas.
%But first we need an elementary lemma to circumvent an extraneous hypothesis in Theorem \ref{Sharifi}. We phrase this lemma in the language of elementary number theory rather than ideals; here ``divides'', ``gcd'' etc. means ``up to associates''.

\begin{Lemma}\label{silly}
Let $U$ be any Unique Factorization Domain and let $d \in U$ with $d \ne 0$.
Suppose $\{ a_m \}_{m=0}^\infty$ is any sequence of nonzero elements of $U$ with $a_m \bigmid a_{m+1}$ for all $m \ge 0$.
Then there exists some $n_0 \ge 0$ such that 
\begin{gather*}
gcd(a_{n_0},d) = gcd(a_m,d) \qquad\text{for all }m \ge n_0, \quad\text{and} \\
gcd(a_m/a_{n_0},d) = 1 \qquad\text{for all }m \ge n_0.
\end{gather*}
\end{Lemma}

\begin{proof}
This is an easy exercise. The key point is that $d$ has only finitely many divisors, so the chain of $gcd(a_m,d)$ must stabilize after finitely many steps.
\end{proof}
\begin{Lemma}
Let $P$ be a finitely generated $\Lambda$-module. 
Let $P_m = \omega_m P$, for all $m \ge 0$. 
Assume there is a $\Lambda$-submodule $N$ of $P$ such that $P_m \subseteq N$ and
$\ord{N/P_m } < \infty$, for all $m \ge 0$.
Assume also that $P/N \iso \ZZ_p$.
Then $P$ is a torsion $\Lambda$-module.
\end{Lemma}

\begin{proof}
By hypothesis $P/N$ is a projective (free) $\ZZ_p$-module, so as $\ZZ_p$-modules we have
\begin{equation}
\label{star}
P/\omega_m P \iso (P/N) \times (N/\omega_m P) \iso \ZZ_p \times finite
\qquad \text{for all } m \ge 0.
\end{equation}
If $P$ is not a torsion $\Lambda$-module, then in Theorem \ref{structure2} we have $r \ge 1$, 
so $P$ has a $\Lambda$-submodule $K$ such that $P/K$ is pseudo isomorphic to $\Lambda$ (where $K$ is the kernel of the map $P\sim \Lambda^r\oplus(\Lambda$-torsion$)\to \Lambda$).
Let overbars denote passage to $P/K$.
Then (as in Claim~1 in Section \ref{4.4} below)
$$
\overline P / \omega_m \overline P \iso \overline {P / \omega_m P} \sim \Lambda / \omega_m \Lambda \iso \ZZ_p[\Gamma_m].
$$
However $\ZZ_p[\Gamma_m]$ is a free $\ZZ_p$-module of rank $p^m$, and so for any $m \ge 1$ it cannot be pseudo-isomorphic to 
a homomorphic image of $P/\omega_m P$ by (\ref{star}) 
and the characterization of finitely generated modules over the PID $\ZZ_p$, a contradiction.

\end{proof}

\begin{Theorem}\label{growth}
Let $P$ be a finitely generated $\Lambda$-module.
Let $P_m = \omega_m P$, for all $m \ge 0$. 
Assume there is a $\Lambda$-submodule $N$ of $P$ such that $P_m \subseteq N$ and $\ord{N/P_m} < \infty$, for all $m \ge 0$. Assume also that $P/N \iso \ZZ_p$.
Then there are nonnegative integers $\mu$ and $\lambda$ and an integer $\nu$ such that 
$$
\ord{N/P_m} = p^{e_m} \qquad\text{where} \quad e_m = \mu p^m + \lambda m + \nu,
$$
for all $m \ge m_0$, for some constant $m_0 \ge 0$.
\end{Theorem}

\begin{proof}
By the preceding lemma, $P$ is a torsion $\Lambda$-module. Let $d = \Char(P)$ and apply Lemma~\ref{silly} in $U = \Lambda$ to $a_m = \omega_m$, for all $m \ge 0$.
Let $n_0$ be as provided by the conclusion of that lemma.
For any $m \ge n_0$, define $\omega_{m,n_0} = \omega_m / \omega_{n_0} \in \Lambda$.
Let $M = P_{n_0}$.\\
\\\
\noindent Note that for all $m \ge n_0$ we have
$$
\omega_{m,n_0} M = (\omega_m / \omega_{n_0}) (\omega_{n_0}P) = \omega_m P = P_m.
$$
By hypotheses then, for all $m \ge n_0$, 
\begin{align*}
\ord{M / \omega_{m,n_0}M} &= \ord{P_{n_0}/P_m}\\
& =\frac{|N/P_m|}{|N/P_{n_0}|} \le \ord {N/P_m} < \infty.
\end{align*}
By Lemma \ref{silly} we have that $\omega_{m,n_0} = \omega_m / \omega_{n_0}$ is relatively prime to $\Char({P})=d$,
for all $m \ge n_0$.
By Proposition~\ref{forras}(3) we have that $\omega_{m,n_0}$ is relatively prime to $\Char(M)$ as well.\\
\\
We now have the hypotheses of Theorem \ref{Sharifi} above. This theorem proves that there are $\mu$, $\lambda$, and some $\nu'$ such that 
$$
\ord{M/\omega_{m,n_0}M} = p^{e'_m} \qquad\text{where} \quad e'_m = \mu p^m + \lambda m + \nu',
$$
for all $m$ greater than or equal to some fixed $m_0 \ge n_0$.\\
\\
Now, as noted above, $\omega_{m,n_0}M = P_m$, and so for all $m \ge m_0$, by Lagrange we have
\begin{align*}
\ord{N/P_m} &= \ord{N/P_{n_0}} \cdot \ord{P_{n_0}/P_m} \\
& = \ord{N/M} \cdot \ord{M/\omega_{m,n_0}M} \\
&= p^k \cdot p^{e'_m} \qquad\text{where} \quad  p^k = \ord{N/M} \quad \text{and}\quad e'_m = \mu p^m + \lambda m + \nu'.
\end{align*}
Finally, let $\nu = k + \nu'$ to obtain the conclusion to the theorem.
\end{proof}

\noindent The goal of the next subsection is to construct a finitely generated torsion $\Lambda$-module, $P=\Pic_\Lambda$. We will then apply Theorem \ref{growth} to $\Pic_\Lambda$ in Section \ref{4.4}.

\subsubsection{Constructing a Finitely Generated $\Lambda$-Module}
\label{4.3.2}

Let $R = \ZZ_p[\Gamma]$ be the usual group ring of $\Gamma$ with coefficients from $\ZZ_p$.
For the given voltage $p$-tower let $X_{p^\infty}$ be its completion, 
so by Theorem~\ref{towercompletion} we may henceforth identify the intermediate graphs
of $X_{p^\infty}/X$ with corresponding graphs in the tower.\\
\\
Fix the following subset of $X_{p^\infty}$:
$$
\mathcal B = \{v_{i,0} \mid 1 \le i \le n \},
$$
where 0 is the additive identity of $\ZZ_p$, so these vertices are taken from the ``zeroth sheet''.
We fix the identification of $X$ and $X_0$ with $\mathcal B$ by $v_i$ is identified with $v_{i,0}$.\\
\\
We first take the free $\ZZ$-module on basis $\mathcal B$, $\Div_{\ZZ}(X)$,
and extend scalars (see \cite{11}, Section~10.4, Corollary~18) to the free $\ZZ_p$-module with the same basis,
now viewed over $\ZZ_p$.
Denote this module by $\Div_{\ZZ_p}(X_0)$.
We can do likewise for each of the graphs $X_m$ and for $X_{p^\infty}$ too.  We obtain the free $\ZZ_p$-modules of divisors 
\begin{align*}
\Div_{\ZZ_p}(X_m) &= \ZZ_p \otimes_{\ZZ} \Div_{\ZZ}(X_m), \quad m \ge 0 \\
\Div_{\ZZ_p}(X_{p^\infty}) &= \ZZ_p \otimes_{\ZZ} \Div_{\ZZ}(X_{p^\infty}) .
\end{align*}

\noindent Now for every $m \ge 0$, each $\Div_{\ZZ}(X_m)$ is a free $\ZZ[\Gamma_m]$-module on the set
$\mathcal B$ too, once we consider the group indices for vertices in $X_m$ to be $p$-adic indices reduced to $\ZZ_p/p^m \ZZ_p \iso \ZZ/p^m \ZZ$; and so $\Div_{\ZZ_p}(X_m)$ is a free module of rank~$n$ over $\ZZ_p[\Gamma_m]$.
We may do likewise for $X_{p^\infty}$ to obtain that $\Div_{\ZZ_p}(X_{p^\infty})$ is a free $R$-module, also of rank~$n$ (on basis $\mathcal B$). In order to emphasize the free, rank~$n$ nature of these respective modules, we adopt the following notation:
$$
\Div_{R_m} = \Div_{\ZZ_p}(X_m) \qquad\text{and}\qquad
\Div_R = \Div_{\ZZ_p} (X_{p^\infty}).
$$

\noindent Since $\Div_{R_m} = \Div_{\ZZ_p}(X_m)$ is a free $\ZZ_p$-module on the basis of vertices of $X_m$, 
$\{v_{i,g} \mid 1 \le i \le n, \ g \in \Gamma_m \}$, 
we may define the usual degree zero divisors with respect to this $\ZZ_p$-basis, and denote this by
$$
\Div^0_{\ZZ_p}(X_m) = \left \{ \sum_{i,g} a_{i,g} v_{i,g} \mid a_{i,g} \in \ZZ_p \text{ and } \sum_{i,g} a_{i,g} = 0 \right \}
$$
where these sums are for $1 \le i \le n$ and $g \in \Gamma_m$.\\
\\
Next we extend scalars from $R$ to $\Lambda$. Since $\Div_R$ is a free $R$-module of rank $n$, its extension is a free
$\Lambda$-module of rank~$n$, denoted by
$$
\Div_\Lambda = \Lambda \otimes_R \Div_R.
$$
Since $R$ is a subring of $\Lambda$ we may simply view the elements of $\Div_\Lambda$ as $\Lambda$-linear combinations
of $\mathcal B$ and $\Div_R$ as the subset of these consisting of $R$-linear combinations of $\mathcal B$.\\
\\
Next we define the Laplacian endomorphism:
$$
\mathcal L_{p^\infty} : \Div_R \longrightarrow \Div_R
\qquad \text{by}\qquad
\mathcal{L}_{p^\infty}(v_{i,0}) = p_{i,0} \quad 1 \le i \le n,
$$
where $p_{i,0}$, the principal divisor ``based at $v_{i,0}$'' is, by definition,
$$
p_{i,0} = n_i \; v_{i,0} - \sum_{\substack{j=1 \\ v_i\sim v_j}}^n v_{j,0+\alpha_{i,j}},
$$
where $n_i$ is the degree of $v_i$ in $X$.
This is extended by $R$-linearity to all of $\Div_R$.
Because $\Gamma$ acts transitively on vertices in each fiber of $X_{p^\infty}/X$, as usual we have that the image of $\mathcal L_{p^\infty}$ 
is the $\ZZ_p$-span of the set of {\it all} principal divisors.  
We encapsulate this by the following notation (definition):
$$
\Pr_R = \mathcal{L}_{p^\infty}(\Div_R).
$$
By taking the ``same map'', but defined on the basis $\mathcal B$ of the free $\Lambda$-module $\Div_\Lambda$ we denote this by 
$$
\widehat{\mathcal L}_{p^\infty} : \Div_\Lambda \longrightarrow \Div_\Lambda
\qquad \text{by}\qquad
\widehat{\mathcal{L}}_{p^\infty}(v_{i,0}) = p_{i,0} \quad 1 \le i \le n,
$$
extended now by $\Lambda$-linearity.  (Formally, $\widehat{\mathcal L}_{p^\infty} = 1 \otimes \mathcal L_{p^\infty}$.)
Now we just define
$$
\Pr_\Lambda = \widehat{\mathcal{L}}_{p^\infty}(\Div_\Lambda).
$$
Likewise, because $\Gamma_m$ acts transitively on the
vertices of $X_m$, using the same $\mathcal L_{p^\infty}$, but instead
reading the vertices $v_{i,0}$ as lying in $\Div_{R_m}$ (i.e., with the vertex indices reduced to
$\ZZ_p / p^m \ZZ_p$), and extended by 
$R_m$-linearity---call this map $\mathcal L_m$---defines the usual Laplacian endomorphism of $\Div_{R_m}$.
Its image is the $R_m$-module of principal divisors of $\Div_{R_m}$, denoted as
$$
\Pr_{R_m} = \mathcal L_m(\Div_{R_m}).
$$
We now define the appropriate Picard groups as follows:
\begin{align*}
\Pic_{R_m}&=\Div_{R_m}/\Pr_{R_m} \qquad &&\text{(an $R_m$-module)}\\
\Pic_R & = \Div_R / \Pr_R \qquad &&\text{(an $R$-module)}\\
\Pic_\Lambda &= \Div_\Lambda / \Pr_\Lambda \qquad &&\text{(a $\Lambda$-module)}.
\end{align*}
So these modules are {\it cokernels} of the respective module endomorphisms.\\
\\
\noindent Next, we identify the $\Lambda$-submodule that plays the role of ``degree zero divisors'' in the proof of Theorem \ref{main}.
\begin{definition}\label{MLambda}
Let 
\begin{align*}
S_1 &=\{v_{i,0}-v_{j,0} \ | \ 1\leq j< i\leq n\} \qquad \text{ and}\\
S_2 &= \{ p_{i,0} \mid 1 \le i \le n \}
\end{align*}
Let $M_\Lambda$ be the $\Lambda$-submodule of $\Div_\Lambda$ generated by $S_1, S_2$ and $(\gamma-1)\Div_\Lambda,$ and let $M_R = \Div_R \cap M_\Lambda$ and $N_\Lambda=M_\Lambda/\Pr_\Lambda.$
\end{definition}

\noindent It turns out that $M_\Lambda$ is actually generated by just $S_1$ and 
$(\gamma-1)\Div_\Lambda$ (see Claim~5(1) in the next subsection).
%but we include the redundant generators for expository clarity.
\noindent Since $\Div_\Lambda$ is a finitely generated $\Lambda$-module and $\Lambda$ is Noetherian, all of its submodules are finitely generated, and so it follows that the quotient modules $\Div_\Lambda/\Pr_\Lambda=\Pic_\Lambda$ and $M_\Lambda/\Pr_\Lambda=N_\Lambda$ are also finitely generated as $\Lambda$-modules.

\begin{Theorem} 
Let $\Div_\Lambda$, $\Pr_\Lambda$, $\Pic_\Lambda$ and $N_\Lambda$ be as above. Then $\Pic_\Lambda$ is a finitely generated module over the Iwasawa Algebra $\Lambda=\ZZ_p[[\Gamma]]$ and therefore so is its submodule $N_\Lambda$.
\end{Theorem}

\subsection{The Main Theorem}
\label{4.4}
We now go on to prove Theorem \ref{main}. 
%The way in which we do so is by relating quotients of $\Lambda$-modules to quotients of $R$-modules. In particular, we show that $M_\Lambda/(\omega_m\Div_\Lambda+\Pr_\Lambda)\iso M_R/(\ker \pi_m +\Pr_R),$ where the latter is isomorphic to $M_{R_m}/\Pr_{R_m}=\mathcal{J}_p(X_m)$, and where $\pi_m:\Div_R\to \Div_{R_m}$. Then by comparing the first column ($\Lambda$-level) with the second column ($R$-level), we are able to show that $N_\Lambda/\omega_m\Pic_\Lambda\iso \mathcal{J}_p(X_m)$. We are then in a position to apply Theorem \ref{growth} which establishes a growth formula for the order of the Jacobians, $\mathcal{J}_p(X_m) \ \forall m.$

\begin{proposition}\label{lattice}  The following diagram holds.

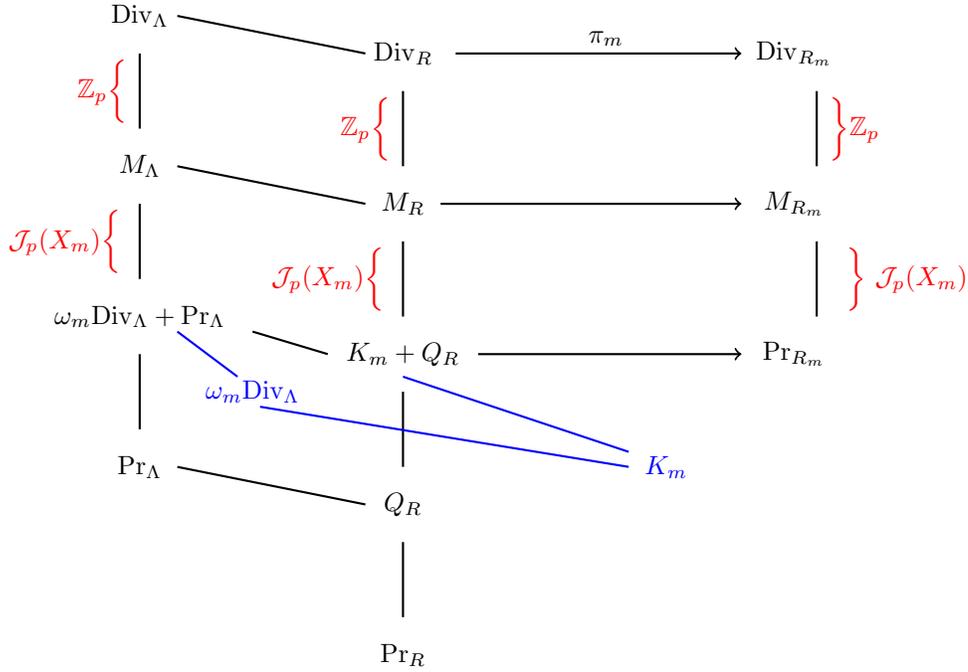
\begin{figure}[H]
\begin{tikzpicture}
%% vertices
%% vertex labels
\node at (-0.5,-4) {$\Pr_R$};
\node at (-0.5,-2) {$Q_R$};
\node at (-4,-1.5) {$\Pr_\Lambda$};
\node at (-0.5,0) {$\textcolor{black}{K_m+Q_R}$};
\node at (-4,.5) {$\textcolor{black}{\omega_m\Div_\Lambda+\Pr_\Lambda}$};
\node at (-2.5,-.5) {$\textcolor{blue}{\omega_m\Div_\Lambda}$};
\node at (3,-1.5) {$\textcolor{blue}{K_m}$};
\node at (4.7,0) {$\Pr_{R_m}$};
\node at (-.5,4) {$\Div_R$};
\node at (-1,3) {$\textcolor{red}{\textcolor{red}{\ZZ_p \bigg\{} }$};
\node at (-4,4.5) {$\Div_\Lambda$};
\node at (4.7,4) {$\Div_{R_m}$};
\node at (5.5,3) {$\textcolor{red}{\bigg\} \ZZ_p}$};
\node at (-.5,2) {$M_R $};
\node at (-1.5,1) {$\textcolor{red}{\mathcal{J}_p(X_m)\bigg\{}$};
\node at (-4,2.5) {$M_\Lambda $};
\node at (-5,1.5) {$\textcolor{red}{\mathcal{J}_p(X_m)\bigg\{} $};
\node at (-4.5,3.5) {$\textcolor{red}{\ZZ_p \bigg\{ }$};
\node at (4.7,2) {$M_{R_m}$};
%\node at (2.2,2.2) {$\pi_m$};
\node at (2.2,4.2) {$\pi_m$};
%\node at (2.2,0.2) {$\pi_m$};
%\node at (2.2,1.2) {$\overline{\pi}_m$};
\node at (6.2,1) {$\textcolor{red}{\bigg\} \ \mathcal{J}_p(X_m)}$};
%%% edges
  %\draw[thick] (-.5,3.5) to (-.5,2.5);
    \draw[thick] (-.5,1.5) to (-.5,.5);
       % \draw[thick] (5,3.5) to (5,2.5);
                \draw[thick] (5,1.5) to (5,.5);
  %\draw[thick][->] (0,0) to (4,0);
    \draw[thick][->] (0,2) to (4,2);
      %\draw[thick][green, ->,->] (0,1) to (4,1);
          \draw[thick] (-.5,2.5) to (-.5,3.5);
           \draw[thick] (5,2.5) to (5,3.5);
               \draw[thick][->] (.2,4) to (4,4);
                  \draw[thick][->] (.5,0) to (4,0);
                    \draw[thick] (-.5,-1.5) to (-.5,-.5);
                      \draw[thick] (-.5,-3.5) to (-.5,-2.5);
                        \draw[thick,blue] (-.5,-.3) to (2.5,-1.3);
                            \draw[thick] (-3.5,4.5) to (-1,4);
                            \draw[thick] (-3.5,2.5) to (-1,2);
                             \draw[thick] (-2.5,.3) to (-1.5,0);
                               \draw[thick] (-3.5,-1.5) to (-1,-2);
                                  \draw[thick](-4,4) to (-4,3);
                                  \draw[thick](-4,2) to (-4,1);
                                  \draw[thick](-4,0) to (-4,-1);
                                   \draw[thick,blue](-3.5,.3) to (-2.7,-.3);
                                      \draw[thick,blue](-2.4,-.7) to (2.5,-1.5);
    
      [decoration={brace,amplitude=7pt}] 
\end{tikzpicture} 
\caption{Lattice and map diagram showing that $M_\Lambda/(\omega_m\Div_\Lambda+\Pr_\Lambda)\iso \mathcal{J}_p(X_m)$,}
\label{figure 5.3}
\end{figure}
\end{proposition}
\noindent Proposition \ref{lattice} is proved by combining the following six claims 1-6 concerning the columns of Figure \ref{figure 5.3}.\\
\\
\noindent First consider the reduction map 
$$
\pi_m: \Div_R\to \Div_{R_m} \quad \text{by} \quad v_{i,g}\mapsto v_{i,\overline{g}}
$$
where $g\in \Gamma$ and $\overline{g}\in \Gamma_m$ is the reduction of $g$ to $\Gamma/\Gamma^{p^m}\cong \ZZ_p/p^m\ZZ_p$, (and recall $\Div_{R_m}= \Div_{\ZZ_p}(X_m)$). Here we are really defining $\pi_m$ on the free $R$-basis vectors on the zeroth sheet, and then extending by $R$-linearity to all of $\Div_R$. It is helpful to keep in mind that for all $m \ge 0$, by the above map and by the previous subsection we have
\begin{gather*}
\text{$\Div_R$ is an $R$-submodule of $\Div_\Lambda$, and} \\
\text{$\Div_{R_m}$ is an $R$-quotient module of $\Div_R$.}
\end{gather*}

\noindent Now let $D = \Div_\Lambda$ as in Proposition \ref{proposition9}, but we simplify notation by writing $\Omega_m(\Div_R)$ to denote $\Omega^D_m(\Div_R).$\\

\noindent \underline{\textbf{Claim 1}}\\
The kernel of $\pi_m$ is $\Omega_m(\Div_R)=\Div_R\cap \omega_m \Div_\Lambda$, where $\omega_m=\gamma^{p^m}-1.$
\begin{proof}
By Proposition \ref{proposition9} and Corollary \ref{corollary6}, we get
the following isomorphisms, where the composition of these isomorphisms is the induced map on $\Div_R$ mod $\ker \pi_m$:
\begin{align*}
\Div_R / \Omega_m(\Div_R) &  \cong   \Div_\Lambda / \omega_m \Div_\Lambda \\
 & \cong (\Lambda\oplus \Lambda\oplus \cdots \oplus \Lambda)/(\omega_m(\Lambda\oplus\Lambda\oplus\cdots\oplus \Lambda))\\
&\cong (\Lambda/(\omega_m))\oplus \cdots \oplus (\Lambda/(\omega_m))\\\
& \cong R_m\oplus \cdots \oplus R_m\\
&\cong \ZZ_p[\Gamma_m]\oplus \cdots\oplus \ZZ_p[\Gamma_m]\\
&\cong \Div_{R_m}
\end{align*}
the free $\ZZ_p[\Gamma_m]$-module of rank $n.$
Thus, the kernel of $\pi_m$ is $\Omega_m(\Div_R).$\\
\end{proof}
\noindent Now let
$$
K_m=\ker \pi_m \qquad \text{and } \qquad Q_R=\Pr_\Lambda \cap \Div_R.
$$

\noindent \textbf{\underline{Claim 2:}}\\
Columns 1 and 2 have the following intersections:\\
\begin{enumerate}
\item $\Div_R\subseteq \Div_\Lambda$\\
\item $\Div_R\cap M_\Lambda=M_R$\\
\item $\Pr_\Lambda\cap \Div_R=Q_R$\\
\item $(\omega_m \Div_\Lambda+\Pr_\Lambda)\cap \Div_R=K_m+Q_R=K_m+\Pr_R$\\

\end{enumerate}
\begin{proof}
(1) holds by Subsection \ref{4.3.2} and (2) and (3) are by definition of $M_R$ and $Q_R$, respectively. 
By Proposition \ref{proposition9} applied with $D = \Div_\Lambda$, $A = \Pr_\Lambda + \omega_m \Div_\Lambda$ and $B = \Pr_R$, since $\Pr_R$ and $\Pr_\Lambda$ are both generated (as $R$- and $\Lambda$-modules, respectively) by the same generators, they both have the same image in $\Div_\Lambda/\omega_m\Div_\Lambda$ as in Claim 1. So, by the {last sentence} of Proposition \ref{proposition9},
\begin{equation}\label{equal}
\Pr_R+\omega_m\Div_\Lambda = \Pr_\Lambda+\omega_m\Div_\Lambda  .
\end{equation}
Then since
$$
\Pr_R\subseteq Q_R \subseteq \Pr_\Lambda,
$$
by (\ref{equal}), we get
\begin{equation}\label{equal1}
\Pr_R+\omega_m\Div_\Lambda =Q_R+\omega_m\Div_\Lambda= \Pr_\Lambda+\omega_m\Div_\Lambda.
\end{equation}
Now because $\Pr_R$ and $Q_R$ are contained in $\Div_R$, intersecting the subgroups in (\ref{equal1}) with $\Div_R$ gives
\begin{align*}
(\Pr_R+\omega_m\Div_\Lambda) \cap \Div_R &= \Pr_R+(\omega_m\Div_\Lambda\cap \Div_R)=\Pr_R+K_m\\
&=Q_R+(\omega_m\Div_\Lambda\cap \Div_R)=Q_R+K_m\\
&=(\Pr_\Lambda+\omega_m\Div_\Lambda) \cap \Div_R,
\end{align*}
which gives (4).
\end{proof}

\noindent \textbf{\underline{Claim 3:}} \\
\noindent {Columns 1 and 2 have the following containments:}\\
\begin{enumerate}
\item $\Pr_\Lambda\subseteq \omega_m\Div_\Lambda+\Pr_\Lambda\subseteq M_\Lambda\subseteq \Div_\Lambda$\\
\item $\Pr_R\subseteq Q_R \subseteq K_m+Q_R \subseteq M_R\subseteq \Div_R$\\

\end{enumerate}

\begin{proof} (1) is clear. 
From Claim 1, we have that $\ker \pi_m =K_m= \Div_R \cap \omega_m \Div_\Lambda$.
Then since $\omega_m = (\gamma-1)(1 + \gamma + \cdots + \gamma^{p^m-1})$, we have
$$
\omega_m \Div_\Lambda = (\gamma-1)(1 + \gamma + \cdots + \gamma^{p^m-1})\Div_\Lambda \subseteq (\gamma-1) \Div_\Lambda \subseteq M_\Lambda.
$$
Thus we have $\ker \pi_m \subseteq \Div_R \cap M_\Lambda = M_R$. So all containments in (2) are clear.
\end{proof}

\noindent \textbf{\underline{Claim 4:}} \\
{Columns 1 and 2 have the following joins:}\\
\begin{enumerate}
\item $M_\Lambda+\Div_R = \Div_\Lambda$
\item $(\omega_m \Div_\Lambda+\Pr_\Lambda)+M_R=M_\Lambda$
\end{enumerate}
\begin{proof}
Apply Proposition \ref{proposition9} to $D = A = \Div_\Lambda$ and $B = \Div_R$.
Its final assertion gives that $\omega_m \Div_\Lambda+\Div_R=\Div_\Lambda.$ Then (1) is immediate because $\omega_m \Div_\Lambda \subseteq M_\Lambda$ by definition. Finally, since $\omega_m \Div_\Lambda+\Div_R=\Div_\Lambda$ and by the latter observation, we have
\begin{align*}
M_\Lambda &= M_\Lambda \cap ( \omega_m \Div_\Lambda + \Div_R) \\
&= \omega_m \Div_\Lambda + (M_\Lambda \cap \Div_R) \\
&= \omega_m \Div_\Lambda + M_R
\end{align*}
as needed for (2). 

\end{proof}

\bigskip\noindent
\noindent \textbf{\underline{Claim 5:}} \\
{As $R$-modules we have the following:}
\begin{enumerate}
\item
$\Div_\Lambda / M_\Lambda  \iso  \ZZ_p$ and $M_\Lambda$ is the $\Lambda$-submodule of $\Div_\Lambda$ generated\\ by $S_1 \cup (\gamma-1)\Div_\Lambda$ (where $S_1$ is as in Definition \ref{MLambda})
\item
$\Div_R / M_R \iso \ZZ_p$
\end{enumerate}

\begin{proof}
To prove (1):
We may obtain $\Div_\Lambda /  M_\Lambda$ as follows: 
first factor   $\Div_\Lambda$   by   $(\gamma - 1)\Div_\Lambda$. By Corollary \ref{corollary6} applied with $m=0$, and as in the proof of Claim 1, we have that
$$
\Div_\Lambda  /  (\gamma - 1) \Div_\Lambda \iso (\Lambda / \omega_0 \Lambda) \oplus \cdots \oplus (\Lambda / \omega_0 \Lambda)   \iso  \underbrace{\ZZ_p \oplus \cdots \oplus \ZZ_p}_{n \text{ of these }},   
$$
where the divisors $v_{1,0}, \dots , v_{n,0}$ map to a basis of this free $\ZZ_p$-module of rank $n$.
Now factor out the submodule generated by the images of all $v_{i,0}-v_{j,0} \ \forall i,j$ from the quotient $\Div_\Lambda/(\gamma-1)\Div_\Lambda$. By doing this, we are simply identifying all the basis vectors with each other,
leaving the rank-1  $\ZZ_p$-module quotient.
The latter process is the same as modding $\Div_{\ZZ_p}(X_0)$ by the degree zero divisors in $\Div_{\ZZ_p}(X_0)$. So this tells us that $p_{i,0}$ for $1\leq i \leq n$ must already be contained in the $\Lambda$-submodule generated by just $S_1\cup (\gamma-1)\Div_\Lambda$.
This proves (1).\\
\\
\noindent To prove (2): First note that by definition $M_R = M_\Lambda \cap \Div_R$. Then by Claim 4(1) we have that 
$$
 \Div_\Lambda / M_\Lambda = (\Div_R + M_\Lambda)/M_\Lambda.
$$
Then by the Diamond Isomorphism Theorem (which says that $(\Div_R+M_\Lambda)/M_\Lambda \iso \Div_R/M_R$) and the previous part we see that 
$$
 \Div_R / M_R\cong \ZZ_p.
$$
as desired.

\end{proof}

\noindent By the proof of Claim 3(2), it follows that the kernel of the map $\pi_m$ restricted to $M_R$ (which we simply denote by $\pi_m$ too) is also equal to $K_m = \Div_R \cap \omega_m \Div_\Lambda$.  From the definition in Section~\ref{4.3.2} we introduce the new notation:
$$
M_{R_m} = \Div^0_{\ZZ_p} (X_m).
$$

\noindent \textbf{\underline{Claim 6:}} \\
\noindent {For Columns 2 and 3 the following holds:}\\
\begin{enumerate}
\item$\Pr_{R_m}\subseteq M_{R_m}\subseteq \Div_{R_m}$\\
\item $\pi_m: \Div_R\to \Div_{R_m}$ is well-defined and surjective\\
\item $\pi_m(M_R)=M_{R_m}$ and $M_R=\pi^{-1}_m(M_{R_m})$\\
\item $\pi_m(K_m+Q_R)=\Pr_{R_m}$ and $\pi^{-1}_m(\Pr_{R_m})=K_m+Q_R$\\
\end{enumerate}

\begin{proof}
(1) and (2) are clear by their respective definitions. 
It is clear that 
$$
S_1 \cup (\gamma-1)\Div_R \subseteq \Div_R \cap M_\Lambda = M_R.
$$
Next we show that the image of $S_1 \cup (\gamma-1)\Div_R$ under $\pi_m$ generates $\Div^0_{\ZZ_p}(X_m)$
as a $\ZZ_p$-module.
The $\ZZ_p$-module $\Div^0_{\ZZ_p}(X_m)$ is generated as a $\ZZ_p$-module by differences of vertices,
$v_{i,r} - v_{j,s}$, where $1 \le i,j \le n$ and $r,s \in \ZZ/p^m \ZZ$. Such differences can be written as 
$$
v_{i,r} - v_{j,s} = (v_{i,r} - v_{i,0})  +  (v_{i,0} - v_{j,0}) +  (v_{j,0} - v_{j,s}) ,
$$ 
where the middle term on the right is in $S_1$.
For each $i$ (and $j$) we may express the first (and third, resp.) differences on the right as telescoping sums of divisors of the form
$v_{i,t+\tau} - v_{i,t}$ for all $t$, 
where $\tau = \pi_m(\gamma)$ is any additive generator for $\ZZ/p^m\ZZ$.
The claim then follows since all of the latter differences are in the image of $(\gamma - 1)\Div_R$ under $\pi_m$.
This argument shows that
\begin{equation}\label{eq2}
M_{R_m} \subseteq \pi_m(M_R).
\end{equation}
To show the reverse containment, let $D = \pi_m^{-1}(M_{R_m})$, (the complete preimage).
By basic properties of homomorphisms (part of the Lattice Isomorphism Theorem) and since $\ker \pi_m \subseteq M_R$ we have:
$$
\pi_m^{-1}(\pi_m(M_R)) = M_R + \ker \pi_m = M_R.
$$
By applying $\pi^{-1}$ to (\ref{eq2}) we get $D \subseteq M_R$.
By the Lattice Isomorphism Theorem we have
that $\pi_m$ induces an isomorphism
$$
\Div_R / D \iso \Div_{R_m} / \pi_m(D) = \Div_{R_m} /  M_{R_m} = \Div_{\ZZ_p}(X_m) / \Div^0_{\ZZ_p}(X_m)  \iso \ZZ_p,
$$
where the latter follows from $\deg:\Div_{\ZZ_p}(X_m)\to \ZZ_p.$\\
\\
Since $D \subseteq M_R$ we get that $\Div_R/M_R$ is a quotient $\ZZ_p$-module of the $\ZZ_p$-module $\Div_R / D$.
By Claim 5(2) we also have $\Div_R/M_R \iso \ZZ_p$. This is illustrated in Figure \ref{figure 5.7}:
\begin{figure}[H]
\centering
\begin{tikzpicture}

\node at (-.5,4.5) {$M_R$};
\node at (1,4.5) {$\textcolor{red}{\ZZ_p}$};
\node at (-.5,6.5) {$\Div_R$};
\node at (-1.2,5.5) {$\textcolor{red}{\ZZ_p  \bigg\{}$};
\node at (-.5,2.5) {$D$};

      \draw[thick] (-.5,6) to (-.5,5);
        \draw[thick] (-.5,4) to (-.5,3);
         %\draw[thick] (1,6.5) to (1,2.5);
           %\draw[thick] (1,6.5) to (0,6.5);
              %  \draw[thick] (1,2.5) to (0,2.5);
       \draw[decoration={calligraphic brace,amplitude=8pt}, decorate, line width=1.25pt]  (.2,6.5) node {} -- (.2,2.5);
\end{tikzpicture}
\caption{$\Div_R/M_R\iso \ZZ_p$ is a quotient $\ZZ_p$-module of the $\ZZ_p$-module $\Div_R/D$}
\label{figure 5.7}
\end{figure}
\noindent However, the only $\ZZ_p$-module quotient of $\ZZ_p$ that is also isomorphic to $\ZZ_p$ is the quotient by the zero submodule (this follows by Lemma \ref{ZZp} in Section \ref{4.3})
i.e., we must have $M_R = D$; and so $\pi_m(M_R) = \pi_m(D) =  M_{R_m}$. This gives (3).\\
\\
\noindent Now $\pi_m$ induces the surjective map
$$
\overline{\pi}_m: M_R/\Pr_R \to M_{R_m}/\Pr_{R_m}
$$
where, by definition, $M_R/\Pr_R=N_R$ and $M_{R_m}/\Pr_{R_m}=\mathcal{J}_p(X_m)$ (the isomorphism $M_{R_m}/\text{Pr}_{R_m}\iso \mathcal{J}_p(X_m)$ is obtained by taking $\otimes\ZZ_p$ to $\Div^0(X_m)/\Pr(X_m)=\mathcal{J}(X_m))$. This is defined by the following commutative diagram in Figure \ref{figure 5.8}:
 \begin{figure}[H]
 \centering
\begin{tikzpicture}
%% vertices
%% vertex labels
\node at (-0.5,2) {$M_R$};
\node at (4,2) {$M_{R_m}$};
\node at (5.5,2.4) {$\small{{\text{proj}}}$};
\node at (4,0) {$M_R/\Pr_R=N_R$};
\node at (9,2) {$M_{R_m}/\Pr_{R_m}=\mathcal{J}_p(X_m)$};
\node at (1.5,2.2) {$\pi_m$};
\node at (1.2,.6) {$\small{{\text{proj}}}$};
\node at (6,1.2) {$\overline{\pi}_m$};
%%% edges
  %\draw[thick] (-.5,3.5) to (-.5,2.5);
    %\draw[thick] (-.5,1.5) to (-.5,.5);
       % \draw[thick] (5,3.5) to (5,2.5);
               %\draw[thick] (5,1.5) to (5,.5);
  %\draw[thick][->] (0,0) to (4,0);
    \draw[thick][->] (0,2) to (3,2);
        \draw[thick][->] (4.5,2) to (6.7,2);
            \draw[thick][->] (0,1.8) to (3.3,.3);
             \draw[thick][->]  (4.5,.2) to (7.7,1.6);
      %\draw[thick][green, ->] (0,1) to (4,1);
    
      [decoration={brace,amplitude=7pt}]
\end{tikzpicture} 
\caption{The map $\overline{\pi}_m:N_R\to \mathcal{J}_p(X_m)$ commutes with the natural projection map}
\label{figure 5.8}
\end{figure}

\noindent To prove (4), first invoke Claim~2(4) to obtain that $K_m + Q_R = K_m + \Pr_R$.
Since $K_m$ is the kernel of $\pi_m$, the subgroups $Q_R$ and $\Pr_R$ have the same image under $\pi_m$.
By definition, $\Pr_R$ is generated as an $R$-module by the principal divisors based at the vertices in the zeroth sheet of
$X_{p^\infty}$; and likewise $\Pr_{R_m}$ is generated as an $R_m$-module by the images of these principal divisors in
$\Div_{R_m}$.  Thus $\pi_m$ maps $\Pr_R$, hence also $K_m + \Pr_R$, surjectively onto $\Pr_{R_m}$.
This gives the first assertion of (4).
Furthermore, since $K_m + \Pr_R$ contains the full kernel of $\pi_m$ and $\Pr_R$ maps onto $\Pr_{R_m}$, the Lattice Isomorphism Theorem immediately gives the second assertion of (4).
\end{proof}
\noindent Claims 1-6 prove Proposition \ref{lattice}.\\
\\
\noindent Now for each subgroup $A$ of $\Div_\Lambda$ let $\widetilde A$ denote the image of $A$ under the natural projection map 
$$
\sim: \Div_\Lambda \longrightarrow \Div_\Lambda / \Pr_\Lambda
$$
(which is both a $\Lambda$- and an $R$-module homomorphism). Since $\sim$ is a $\Lambda$-module homomorphism, the image of $\omega_m \Div_\Lambda + \Pr_\Lambda$ 
under it is $\omega_m \Pic_\Lambda$. 
Since $\Div_R$ is an $R$-submodule of $\Div_\Lambda$, we may apply $\sim$ to it as well, and to its submodules.\\
\\
By the Diamond Isomorphism Theorem, since we've checked all the appropriate intersections from column~1 to column~2 in Figure~\ref{figure 5.3}, this natural projection gives the first two columns in Figure~\ref{figure 5.4} as well as all intersections (depicted, as usual, by horizontal lines) between their subgroups in column~2.  To get the third column of Figure~\ref{figure 5.4}, factor the third column of Figure~\ref{figure 5.3} by $\Pr_{R_m}$. The horizontal lines---which are homomorphisms---relating column~2 to column~3 in Figure~\ref{figure 5.4} are obtained by taking images of the subgroups in column~2 under $\overline \pi_m$.  
%By Claim~2(4), $\overline \pi_m$ is a well-defined group homomorphism from the second column of Figure~\ref{figure 5.4} to its third column. By simple inspection of the claims, all the horizontal group homomorphisms from column~2 to their images in column~3 of Figure~\ref{figure 5.4} are valid too.  Note that there are no direct ``horizontal line'' relationships from column~1 to column~3, and so nothing to check in that regard. 
By the Lattice Isomorphism Theorem, the (already established) quotient groups (in red) are consequently also preserved when passing between columns (thus also transitively from column~1 to column~3).  We only need these to be abelian group isomorphisms; but they are, in fact, $R$- and $\Lambda$-module isomorphisms.

\begin{figure}[H]
\begin{tikzpicture}
%% vertices
%% vertex labels
\node at (-0.5,-2) {$\widetilde{0}$};
\node at (-4,-1.5) {$0$};
\node at (-0.5,0) {$\widetilde{K}_m$};
\node at (-4.3,.5) {${\omega_m\Pic_\Lambda}$};
%\node at (-2.5,-.5) {$\textcolor{blue}{\omega_m\Div_\Lambda}$};
%\node at (3,-1.5) {$\textcolor{blue}{K_m}$};
\node at (5,0) {$\overline{0}=\overline{\Pr}_{R_m}$};
\node at (-.5,4) {$\widetilde{\Div}_R$};
\node at (-1,3) {$\textcolor{red}{\textcolor{red}{\ZZ_p \bigg\{} }$};
\node at (-4,4.5) {$\Pic_\Lambda$};
\node at (5,4) {$\Pic_{R_m}$};
\node at (5.5,3) {$\textcolor{red}{\bigg\} \ZZ_p}$};
\node at (-.5,2) {$\widetilde{M}_R $};
\node at (-1.5,1) {$\textcolor{red}{\mathcal{J}_p(X_m)\bigg\{}$};
\node at (-4,2.5) {$N_\Lambda $};
\node at (-5,1.5) {$\textcolor{red}{\mathcal{J}_p(X_m)\bigg\{} $};
\node at (-4.5,3.5) {$\textcolor{red}{\ZZ_p \bigg\{ }$};
\node at (5,2) {$\mathcal{J}_p(X_m)$};
\node at (6,1) {$\textcolor{red}{\bigg\}\mathcal{J}_p(X_m)}$};
%\node at (2.2,2.2) {$\pi_m$};
\node at (2.2,4.3) {$\overline{\pi}_m$};
%\node at (2.2,0.2) {$\pi_m$};
%\node at (2.2,1.2) {$\overline{\pi}_m$};
%\node at (6.2,1) {$\textcolor{red}{\bigg\} \ \mathcal{J}_p(X_m)}$};
%%% edges
  %\draw[thick] (-.5,3.5) to (-.5,2.5);
    \draw[thick] (-.5,1.5) to (-.5,.5);
       % \draw[thick] (5,3.5) to (5,2.5);
                \draw[thick] (5,1.5) to (5,.5);
  %\draw[thick][->] (0,0) to (4,0);
    \draw[thick][->] (0,2) to (4,2);
      %\draw[thick][green, ->,->] (0,1) to (4,1);
          \draw[thick] (-.5,2.5) to (-.5,3.5);
           \draw[thick] (5,2.5) to (5,3.5);
               \draw[thick][->] (.5,4) to (4,4);
                  \draw[thick][->] (.5,0) to (4,0);
                    \draw[thick] (-.5,-1.5) to (-.5,-.5);
                     % \draw[thick] (-.5,-3.5) to (-.5,-2.5);
                        %\draw[thick] (-.5,-.3) to (2.5,-1.3);
                            \draw[thick] (-3.5,4.5) to (-1.5,4);
                            \draw[thick] (-3.5,2.5) to (-1,2);
                             \draw[thick] (-3.5,.3) to (-1,0);
                               \draw[thick] (-3.5,-1.5) to (-1,-2);
                                  \draw[thick](-4,4) to (-4,3);
                                  \draw[thick](-4,2) to (-4,1);
                                  \draw[thick](-4,0) to (-4,-1);
                                  % \draw[thick](-3.5,.3) to (-2.7,-.3);
                                     % \draw[thick](-2.4,-.7) to (2.5,-1.5);
    
      [decoration={brace,amplitude=7pt}]
\end{tikzpicture} 
\caption{The natural projection homomorphism from $\Div_\Lambda$ to $\Pic_\Lambda$ indicated in the first two columns and the passage from $\pi_m$ to $\overline{\pi}_m$ indicated in the third column}
\label{figure 5.4}
\end{figure}
\noindent Now we apply Theorem \ref{growth} with $P=\Pic_\Lambda$ and $N=N_\Lambda= M_\Lambda/\Pr_\Lambda$ used as $P$ and $N$ in its hypotheses. Since $X_m$ is connected and $\omega_m\Pic_\Lambda\subseteq N_\Lambda$ for all $m$ by Figure $\ref{figure 5.4}$, we have
$$
|N_\Lambda/ \omega_m \Pic_\Lambda | = |\mathcal{J}_p(X_m)| <\infty.
$$
%It is here that we need the crucial hypothesis that all $X_m$ are connected, so the Jacobians are finite. 
This leads immediately to the conclusion of Theorem \ref{main}.
\\
\\
An important invariant of any voltage graph with abelian voltage group is \textit{the reduced Stickelberger Element.}  
With respect to the basis $\mathcal B$ of both $\Div_R$ (as an $R$-basis) and 
$\Div_\Lambda$ (as a $\Lambda$-basis) (see Subsection \ref{4.3.2}), the two maps, $\mathcal{L}_{p^\infty}$ and $\widehat{\mathcal L}_{p^\infty}$, have {\it the same matrix representation}.
Thus we define
$$
\Theta_{p^\infty} = \det \mathcal{L}_{p^\infty} = \det \widehat{\mathcal{L}}_{p^\infty},
$$
which is an element of $R = \ZZ_p[\Gamma]$ that plays the role of the \textit{reduced Stickelberger element}.\\
\\
\textit{Remark:} $\Theta_{p^\infty}$ annihilates both $\Pic_R$ and $\Pic_\Lambda$ --- see \cite{29} for discussion and additional uses of the reduced Stickelberger element.

\begin{Corollary}\label{sticky}
Under the hypothesis and notation of Theorem \ref{main}, 
the ranks of $\mathcal J_p(X_m)$ are bounded as $m \rightarrow \infty$
if and only if $p$ does not divide $\Theta_{p^\infty}$ in $\Lambda$ (or in $\ZZ_p[\Gamma]$). 
\end{Corollary}

\begin{proof}
By definition, $\Pic_\Lambda$ is the cokernel of the voltage Laplacian, $\mathcal L_{p^\infty}:\Div_\Lambda\to \Div_\Lambda,$
where $\Theta_{p^\infty} = \det \mathcal L_{p^\infty}$. In the notation of Theorem \ref{structure2}, let $p^\mu$ be the product of the $p^{k_i}$. Then the characteristic polynomial, as in Definition \ref{characteristic}, is equal to
\begin{equation}\label{char}
p^\mu \prod_{j=1}^t g_j^{m_j} = \Char(\Pic_\Lambda).
\end{equation}
Let $M = \omega_{m_0}  \Pic_\Lambda$ where $m_0 \ge 0$ is fixed. Then since $\Pic_\Lambda /M$ has finite $p$-rank (fixed, independent of $m \rightarrow \infty$), the ranks of $\Pic_\Lambda$ and $M$ differ by a constant,
and one is bounded as $m \rightarrow \infty$ if and only if the other is bounded.\\
\\
We now compare $\mu$ invariants for $\Pic_\Lambda$ and $M$, as follows. 
By Proposition~\ref{forras}(3), we have
\begin{equation}\label{char}
\Char(M) = \frac{\Char(\Pic_\Lambda)}{\Char(\Pic_\Lambda/M)} .
\end{equation}
The $\Lambda$-module $\Pic_\Lambda/M$ is a quotient of the module
$\Div_\Lambda / (\omega_{m_0} \Div_\Lambda)$;
and as in Claim~1,
$$
\Div_\Lambda / (\omega_{m_0} \Div_\Lambda) \cong \underbrace{(\Lambda/\omega_{m_0}\Lambda)\oplus\cdots \oplus (\Lambda/\omega_{m_0}\Lambda)}_{n \text{ of these}}.
 $$
But by Lemma \ref{lemma1} we know $\omega_{m_0}$ maps to a distinguished polynomial in $\ZZ_p[[T]] \iso \Lambda$, 
so $(\Lambda/(\omega_{m_0}))^n$ is already in Iwasawa decomposition form, and it clearly has characteristic polynomial 
$\omega_{m_0}^n$ (again, under the identification $\gamma \mapsto T+1$).
One more usage of Proposition \ref{forras}(3) gives that
$$
\Char(\Pic_\Lambda/M) \bigmid \omega_{m_0}^n ,
$$
so $\Char(\Pic_\Lambda/M)$ is relatively prime to $p$ (since the distinguished polynomial $\omega_{m_0}$ is). \\
By (\ref{char}), this shows
$$
p \bigmid \Char(\Pic_\Lambda) \quad \Longleftrightarrow \quad p \bigmid \Char (M).
$$
If $\mu(\Pic_\Lambda) = 0$, then $\Char(\Pic_\Lambda)=\Theta_{p^\infty}$ by Proposition 10.23 in \cite{21}.
In this case, $p$ does not divide $\Theta_{p^\infty}$ by definition of $\Char(\Pic_\Lambda)$. By Lemma 13.20 of \cite{28}, we have that the ranks of the finite $\Lambda$-module quotients of a finitely generated torsion $\Lambda$-module
stay bounded if and only if the $\mu$ invariant of the Iwasawa decomposition is zero.
So if the ranks of $\mathcal J_p(X_m)$ stay bounded as $m \to \infty$, then $p$ does not divide $\Theta_{p^\infty}$.\\
\\
Conversely, we show that if the ranks of $\mathcal J_p(X_m)$ don't stay bounded as $m \rightarrow \infty$, then $p$ does divide $\Theta_{p^\infty}$ in $\Lambda$. So if the rank of $\mathcal{J}_p(X_m) \to \infty$ as $m \to \infty,$
then the $\mu$-invariant of the submodule $M$, and hence also of $\Pic_\Lambda,$ must be nonzero.
i.e., the Iwasawa decomposition of $\Pic_\Lambda$ must have at least one factor of the form
$\Lambda/(p^a)$, for some $a \ge 1.$
This forces $p$ to divide $\Theta_{p^\infty}$ as follows.
By Corollary 9 in \cite {29} or by \cite{28} page 297, $\Theta_{p^\infty}$ annihilates $\Pic_\Lambda$.
It follows from the definition of pseudo-isomorphism that there is some submodule of fixed finite index in the 
Iwasawa Decomposition factor $\Lambda/(p^a)$ that is also annihilated by $\Theta_{p^\infty}$;
and hence $\Theta_{p^\infty}$ annihilates a submodule of fixed finite index in every quotient of $\Lambda/(p^a)$.
But the latter module has quotient modules, $\Lambda/(p^a,T^k)$, of order $p^{ak}$, for every $k\ge 1$, and none of these 
possess a nonzero submodule annihilated by an element of $\Lambda$ that is prime to $p$. 
Thus $p$ must divide $\Theta_{p^\infty}$.
%This completes the proof of the converse.
\end{proof}

\begin{exmp}
Let
$
X=X_0\leftarrow X_1\leftarrow X_2 \leftarrow \cdots \leftarrow X_m \leftarrow \cdots
$
be a cyclic voltage $p$-tower of derived graphs over a base graph $X$ whose vertices are $v_1,\dots,v_n$, 
with $v_1$ adjacent to $v_2$.
Assume that the voltage adjacency matrix for $X_m/X$ (with the standard orientation on $X$) has a generator for the cyclic group
of order $p^m$ in entry $1,2$; its inverse in entry 2,1; and has a $1$ (the identity of the voltage group) in entry
$i,j$ whenever $v_i$ is adjacent to $v_j$ for $\{i,j\} \ne \{1,2\}$; and has zeros elsewhere (including on the diagonal). 
Then 
$$
|\mathcal{J}_p (X_m)| = p^{e_m}, \quad \text{where} \quad e_m=\mu p^m+\lambda m+\nu
\qquad \text{for all } m \ge 0,
$$
where $p^\mu$ is the largest power of $p$ dividing the reduced Stickelberger element 
($p^\mu$ can be shown to be independent of $m$),
and $\lambda = 1$.
In particular, when $X$ is the complete graph on $n$ vertices, $p^\mu$ is the
largest power of $p$ dividing $(n - 2) n^{n-3}$, for $n \ge 3$.

Starting just with a base graph $X$ as above, we can choose generators for voltage groups $\ZZ / p^m \ZZ$
so that the sequence of $1,2$ entries in such voltage adjacency matrices converges in $\ZZ_p$.
This constructs derived graphs $X_m/X$ that form a cyclic voltage $p$-tower (and if $X_1$ is connected, then all $X_m$ are too) --- 
see \cite{29} for further details on these ``single voltage'' derived covers. 
\end{exmp}

\end{document}